\documentclass[journal]{IEEEtran}
\usepackage{amsmath}
\usepackage{amssymb}
\usepackage[hidelinks]{hyperref}
\usepackage{subfigure}
\usepackage{accents}
\usepackage{tcolorbox}
\usepackage{tikz,lipsum}
\tcbuselibrary{skins,breakable}
\usetikzlibrary{calc,arrows}

\usepackage{bm}
\usepackage{mathtools}

\newtheorem{thm}{Theorem}[section]
\newtheorem{corollary}[thm]{Corollary}
\newtheorem{lemma}[thm]{Lemma}
\newtheorem{prop}[thm]{Proposition}

\newtheorem{assumption}[thm]{Assumption}

\newcommand{\real}{\ensuremath{\mathbb{R}}}

\newcommand{\complex}{\ensuremath{\mathbb{C}}}

\newcommand{\subscr}[2]{#1_{\textup{#2}}}
\newcommand{\supscr}[2]{#1^{\textup{#2}}}

\newcommand{\diag}[1]{\operatorname{diag}(#1)}

\newcommand{\dist}{\operatorname{dist}}

\newcommand{\eps}{\varepsilon}
\renewcommand{\epsilon}{\varepsilon}

\newcommand{\Pc}{{\mathcal{P}}}

\usepackage{psfrag}

\newcommand{\RgeO}{\ensuremath{\mathbb{R}_{ \geq 0}}}
\newcommand{\Rat}[1]{\ensuremath{\mathbb{R}^{#1}}}
\newcommand{\B}{\mathcal{B}}

\newcommand{\D}{\mathcal{D}}

\newcommand{\F}{\mathcal{F}}
\renewcommand{\H}{\mathcal{H}}

\newcommand{\N}{\mathcal{N}}
\newcommand{\Obs}{\bm{\mathcal{O}}}
\renewcommand{\P}{\mathcal{P}}
\newcommand{\W}{\mathcal{W}}
\newcommand{\X}{\mathcal{X}}
\newcommand{\bP}{\mathbb{P}}
\newcommand{\bN}{\mathbb{N}}
\newcommand{\fn}{\mathfrak n}
\newcommand{\fm}{\mathfrak m}
\newcommand{\fK}{\mathfrak K}

\newcommand{\bE}{\ensuremath{\mathbb{E}}}

\newcommand{\diameter}[1]{\operatorname{diam}(#1)}
\newcommand{\supp}[1]{{\rm supp}(#1)}

\newcommand{\oprocendsymbol}{\hbox{$\square$}}
\newcommand{\oprocend}{\relax\ifmmode\else\unskip\hfill\fi\oprocendsymbol}

\renewcommand{\baselinestretch}{.99}

\allowdisplaybreaks

\newcommand{\longthmtitle}[1]{\mbox{}{\bf \textit{(#1).}}}

\begin{document}



\title{Data-driven 
  ambiguity sets with probabilistic guarantees for dynamic
  processes\thanks{A preliminary version of this work appeared
    as~\cite{DB-JC-SM:19-ecc} in the European Control Conference. This
    work was supported by the DARPA Lagrange program through award
    N66001-18-2-4027.}}

\author{Dimitris Boskos\qquad Jorge Cort{\'e}s \qquad Sonia
  Mart{\'i}nez\thanks{The authors are with the Department of
    Mechanical and Aerospace Engineering, University of California,
    San Diego, \tt{\{dboskos,cortes,soniamd\}@ucsd.edu}.}}

\maketitle

\begin{abstract}
  Distributional ambiguity sets provide quantifiable ways to
  characterize the uncertainty about the true probability distribution
  of random variables of interest.  This makes them a key element in
  data-driven robust optimization by exploiting high-confidence
  guarantees to hedge against uncertainty.  This paper explores the
  construction of Wasserstein ambiguity sets in dynamic scenarios
  where data is collected progressively and may only reveal partial
  information about the unknown random variable.  For random variables
  evolving according to known dynamics, we leverage assimilated
  samples to make inferences about their unknown distribution at the
  end of the sampling horizon.  Under exact knowledge of the flow map,
  we provide sufficient conditions that relate the growth of the
  trajectories with the sampling rate to establish a reduction of the
  ambiguity set size as the horizon increases.  Further, we
  characterize the exploitable sample history that results in a
  guaranteed reduction of ambiguity sets under errors in the
  computation of the flow and when the dynamics is subject to bounded
  unknown disturbances.  Our treatment deals with both full- and
  partial-state measurements and, in the latter case, exploits the
  sampled-data observability properties of linear time-varying systems
  under irregular sampling.  Simulations on a UAV detection
  application show the superior performance resulting from the
  proposed dynamic ambiguity sets.
\end{abstract}


\IEEEpeerreviewmaketitle
	
\section{Introduction}

In stochastic optimization, ambiguity sets play a key role by
accounting for 'what-if' scenarios regarding the true probability
distribution of random variables affecting the objective function or
the constraints.  Rigorous guarantees on the probability of these sets
containing the true distribution allows the designer to robustify
decisions in the face of uncertainty. This explains the numerous
applications that distributionally robust optimization with ambiguity
sets finds in decision making under uncertainty, reliability-based
design, and data-driven modeling.  Moving beyond the static problem,
where full-state measurements on the random variable are available all
at once, this paper instead looks at scenarios where the random
variable evolves dynamically and data is collected incrementally.  We
are interested in developing methods to construct ambiguity sets and
track their evolution across time while maintaining their
probabilistic guarantees about the true distribution. These methods
should handle exact and approximate knowledge of the dynamics, the
presence of disturbances, and the availability of partial-state
measurements of the random variable. We also seek to shed light on the
trade-offs between data assimilation and accuracy of the resulting
ambiguity sets.

\textit{Literature review:} DRO optimization is an area of stochastic
programming~\cite{AS-DD-AR:14} which has gained significant recent
research attention~\cite{PME-DK:17},~\cite{AS:17}, \cite{DB-MS-MZ:18},
in view of the progress on robust optimization during the last two
decades~\cite{AB-LEG-AN:09}. A main characteristic of DRO is that
worst-case decisions against model uncertainty can be quantified with
performance guarantees, by considering a set of distributions up to a
certain distance from a candidate model. There is an exhaustive number
of choices for distances in spaces of probability distributions
\cite{STR-LK-SVS-FF:13}. Among the most popular distance-type notions
for DRO problems are $\phi$-divergences \cite{AB-DD-AD-BM-GR:13},
\cite{RJ-YG:16}, and Wasserstein metrics \cite{RG-AJK:16},
\cite{JB-KM:19}. For data-driven problems where robustness is measured
with respect to the empirical distribution, the Wasserstein distance
becomes a suitable choice, since it does not require any absolute
continuity conditions between the associated distributions. The
work~\cite{PME-DK:17} leverages recent concentration of measure
inequalities~\cite{NF-AG:15} to build Wasserstein ambiguity sets
around the empirical distribution of the data, and provides tractable
reformulations of the associated DRO problems with out-of-sample
guarantees. These are exploited in~\cite{AC-JC:17-tac}, where a
distributed reformulation of the min-max DRO problem is established
via saddle-point dynamics, and in \cite{DL-SM:18-extended}, where
online sample assimilation is fused with an efficient optimization
algorithm to provide on-the-fly data-driven DRO solutions. It is
worthwhile mentioning also the work \cite{JB-YK-KM:16}, where the
notion of a robust Wasserstein profile function is employed, providing
fast asymptotic convergence rates for high-dimensional samples. Recent
work has considered distributionally robust Kalman filtering
approaches for the state estimation of uncertain time-varying
processes for the Kullback-Leibler~\cite{BCL-RN:13},
$\tau$-divergences~\cite{MZ:17}, and the
Wasserstein~\cite{SSA-VAN-DK-PME:18} metrics.

Observability of linear and nonlinear systems occupies a central part
of the control literature~\cite{EDS:98}. Of considerable practical
interest is the case where the output of a system is not continuously
measured and samples are collected instead. Classical results
regarding observability for linear time-invariant systems under
periodic sampling with equidistributed measurements can be found
in~\cite{EDS:98}. For the same system class, a periodic sampling
schedule which always maintains observability of the continuous plant
was proposed in~\cite{GK:99}. Observability under regular sampling for
nonlinear models has been studied in~\cite{SA-JCV:04} for systems on
compact manifolds, in~\cite{SA-HF-JCV:14} for bilinear systems, and
conditions under which the property become generic are given
in~\cite{SA-MM-JCV:18}. Results on the asymptotic state estimation in
the nonlinear case through sampled-data observes are derived
in~\cite{IK-CK:09}. Also, the observability of linear time-invariant
(LTI) systems under irregular sampling, using properties of
exponential polynomials, has been studied in the recent
works~\cite{LYW-CL-GY-LG-CZX:11,SZ-FA:16}, and exploited to establish
observability of LTI ensembles in~\cite{SZ-HI-FA:17}.
	
\textit{Statement of contributions:} We consider dynamically evolving
random variables of unknown distribution on which samples from
multiple independent realizations are collected in an online fashion.
Our contributions revolve around providing solutions to integrate the
data to construct distributional ambiguity sets based on the
Wasserstein metric that enjoy rigorous probabilistic guarantees.
Throughout the technical approach, we pay attention to characterizing
how the information contained in the incrementally collected data can
be pushed forward in time to infer properties about the evolving
distribution of the process.  Under full-state measurements and exact
pushforward through the flow, our first contribution builds ambiguity
balls that incorporate past data to enjoy desirable guarantees on the
probability of containing the true distribution. We also identify
conditions on the growth of the dynamics of the random variable under
which the ambiguity radius shrinks as the horizon increases.  Our
second contribution considers scenarios with approximate pushforwards
and disturbances in the dynamics, and characterizes the modifications
necessary in the ambiguity radius that retain the same high-confidence
probabilistic guarantees.  This result allows us to quantify the
effective sampling horizon that ensures the monotonic reduction of the
ambiguity set with the number of sampling times.  To enable the
extension of these results to the case with partial-state
measurements, our third contribution studies robust sample-data
observability under irregular sampling of linear time-varying
systems. We provide conditions on the inter-sampling time of the
trajectories under which full-state information can be extracted. Our
final contribution is the construction of high-confidence
distributional ambiguity sets under partial-state measurements. We
illustrate our results in a UAV detection scenario application
formulated as a distributionally robust optimization
problem\footnote{Throughout the paper, we use the following
  notation. We denote by $\|\cdot\|$ and $\|\cdot\|_{\infty}$ the
  Euclidean and infinity norm in $\Rat{n}$, respec., and by
  $[n_1:n_2]$ the set of integers
  $\{n_1,n_1+1,\ldots,n_2\}\subset\bN\cup\{0\}$.  The interpretation
  of a vector in $\Rat{n}$ as an $n\times 1$ matrix should be clear
  form the context (this avoids writing double transposes).  Given a
  differentiable function $G:\Rat{n}\to\Rat{m}$, $DG(x)$ denotes its
  derivative at $x\in\Rat{n}$. The diameter of $S\subset\Rat{n}$ is
  $\diameter{S}:=\sup\{\|x-y\|_{\infty}\,|\,x,y\in S\}$ and the
  distance of $x\in\Rat{n}$ to $S$ is
  $\dist(x,S):=\inf\{\|x-y\|\,|\,y\in S\}$.  For $A\in\Rat{m\times
    n}$, $A^{\dagger}$ denotes its Moore-Penrose pseudoinverse, and
  $\|A\|$ its induced Euclidean norm, namely,
  $\|A\|:=\max_{\|x\|=1}\|Ax\|/\|x\|$. For the distinct eigenvalues
  $\lambda_1,\ldots,\lambda_q$ of $A$, the index $\fm_i$ of each
  eigenvalue is the exponent of the term $\lambda-\lambda_i$ in the
  minimal polynomial
$\prod_{j=1}^q(\lambda-\lambda_j)^{\fm_j}$ of $A$;
equivalently, $\fm_j$ is the dimension of the largest Jordan block
corresponding to $\lambda_j$.
We let $A\otimes B$ denote the Kronecker product. We denote by
$\diag{a_1,\ldots,a_n}$ the diagonal matrix with entries
$a_1,\ldots,a_n$ in its main diagonal. For any $z\in\complex$,
$\Im(z)$ represents its imaginary part.
We denote by $\B(\Rat{d})$ the Borel $\sigma$-algebra on $\Rat{d}$,
and by $\P(\Rat{d})$ the space of probability measures on
$(\Rat{d},\B(\Rat{d}))$.  Given a real number $p\ge 1$, we denote by
$\P_p(\Rat{d})$ the set of probability measures in $\P(\Rat{d})$ with
finite $p$-th moment, i.e.,
$\P_p(\Rat{d}):=\{\mu\in\P(\Rat{d})\,|\,\int_{\Rat{d}}\|x\|^pd\mu<\infty\}$.
For any $p\ge 1$, and probability measures $\mu,\nu\in\P_p(\real^d)$,
the Wasserstein distance is
\begin{align*}
  W_p(\mu,\nu):=\Big(\inf_{\pi\in\H(\mu,\nu)} \Big\{
  \int_{\Rat{d}\times\Rat{d}}\|x-y\|^p \pi(dx,dy)\Big\}\Big)^{1/p},
\end{align*}
where $\H(\mu,\nu)$ is the set of all probability measures on
$\Rat{d}\times\Rat{d}$ with marginals $\mu$ and $\nu$,
respectively. For any $\mu\in\P(\Rat{d})$, its support is the closed
set $\supp{\mu}:=\{x\in\Rat{d}\,|\,\mu(U)>0\;\textup{for each
neighborhood}\;U\;{\rm of}\;x\}$,
or equivalently, the smallest closed set with measure one.  Given two
measurable spaces $(\Omega,\F)$, $(\Omega',\F')$, and a measurable
function $\Psi$ from $(\Omega,\F)$ to $(\Omega',\F')$, the pushforward
map $\Psi_{\#}$ assigns to each measure $\mu$ in $(\Omega,\F)$, a new
measure $\nu$ in $(\Omega',\F')$, defined by $\nu:=\Psi_{\#}\mu$ iff
$\nu(B)=\mu(\Psi^{-1}(B))$ for all $B\in\F'$. The map $\Psi_{\#}$ is
linear and satisfies $\Psi_{\#}\delta_{\omega}=\delta_{\Psi(\omega)}$,
with $\delta_{\omega}$ the Dirac measure centered at
$\omega\in\Omega$.}.

\section{Problem Formulation}\label{sec:pf}
	
A distributionally robust optimization problem (DROP) takes the
form
\begin{align*} 
{\rm (P)}\qquad\qquad\inf_{x\in\X}\sup_{P\in\widehat{\P}^N}\bE_P[f(x,\xi)],
\end{align*}
where $x\in\X\subset\Rat{n}$ is the decision variable and $\xi$
represents a random variable distributed according to a distribution
$P_\xi\in \Pc(\real^d)$.  This distribution is unknown and hence one
formulates a worst-case expectation problem over an \textit{ambiguity
  set}~$\widehat{\P}^N$ which contains it with some probabilistic
guarantee. Such guarantees can be obtained when data about the random
variable is available: given $N$ i.i.d.~samples $\xi^1,\ldots,\xi^N$
drawn according to the unknown distribution $P_\xi$ and a reliability
parameter $\beta \in (0,1)$, one can construct~\cite{PME-DK:17} an
ambiguity set
such that $\bP(P_{\xi}\in\widehat{\P}^N)\ge 1-\beta$.  We refer to
this common formulation as a \textit{static} DROP.  Here, we are
interested in building the ambiguity set from dynamically varying,
possibly partial data in an online manner, as it may not be feasible
to collect and process many samples in a given time instant. We next
present an illustrative example.
	
\subsection{Motivating example}\label{subsec:motivating:example}
	
\begin{figure}[tbh]
  \centering
  \includegraphics[width=.9\columnwidth]{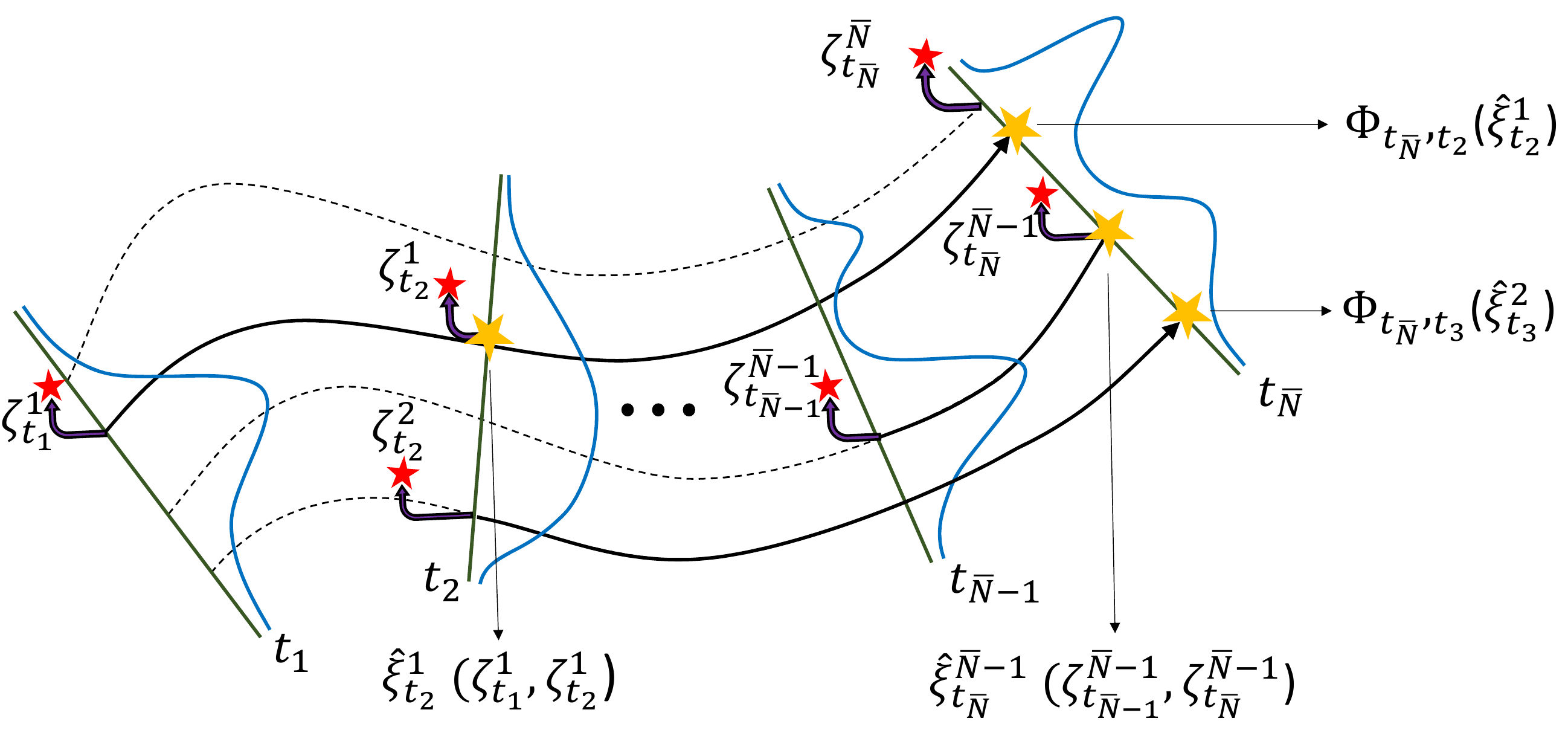}	
  \caption{Partial measurements across a sequence of time instants of
    trajectories of the random variable can be exploited to generate
    samples of the full-state distribution at $T$. The latter can be
    employed to construct ambiguity sets that contain the true unknown
    distribution with high confidence.  The blue curves show how the
    probability density of the state distribution evolves over
    time. We display the output samples with the red stars and use a
    bent arrow to represent the output map applied to the
    trajectories' states at the sampling instants. Each yellow star
    depicts either the reconstructed state of a trajectory once the
    last output sample is collected from it, or its corresponding
    value pushed forward to~$t_{\bar N}=T$.
  }\label{fig:output:samples:from:trajectories}
\end{figure}

Let $\xi_t=(\xi_{1t},\xi_{2t})$ describe the position and velocity of
a unit acceleration particle, evolving according to the known dynamics
\begin{align*}
  \dot{\xi}_{1t}=\xi_{2t},\quad \dot{\xi}_{2t} =1,
\end{align*}
over a time horizon $[0,T]$. Assume that we can measure the position
$\zeta_t=H(\xi_t):=\xi_{1t}$, but not the velocity. Then, taking two position 
measurements are sufficient to reconstruct the full state of the particle. Both 
its initial position and velocity are random with unknown
probability distribution.  Therefore, the particle state at each $t$
is a random variable with law $P_{\xi_t}$.
Given that these distributions are unknown, we focus on specifying an
ambiguity set $\widehat{\P}_T^N$ at time $T$ which contains the true
distribution $P_{\xi_T}$ with high confidence. To do this, independent
output samples from the distribution of $\zeta_t=H(\xi_t)$ are
available at time instants $0\le t_1<\ldots<t_{\bar N}=T$, cf.
Figure~\ref{fig:output:samples:from:trajectories}.  In this scenario,
a direct application of the conventional approach to static DROPs
would only employ the  samples obtained at $t_{\bar N}$
to construct the ambiguity set and solve (P).
This would not be desirable because (i) constructing reliable
ambiguity sets requires a finite, but sufficiently large amount of
data, and (ii) the output map reveals information about the particle's
position, but not about its velocity.

%

These considerations motivate the question of how to leverage the full
set of samples and the system's observability properties to construct
a better ambiguity set at time~$T$.  If full-state samples were drawn
from different realizations of the system at each $t_i$, they could be
pushed forward through the flow map to obtain their corresponding
values at $T$ and increase the exploitable data for the construction
of the ambiguity set. Building on this observation, we address here
three interconnected problems: (i) the dynamics transforms, and may
potentially increase, the initial uncertainty about the random
variable. It is therefore of interest to understand to what extent
this may be compensated by the number of collected samples; (ii)
bounded errors in the dynamics or, even if it is fully known, when the
flow map cannot be computed exactly, induce errors in the propagation
of the collected samples. It is therefore of interest to characterize
accuracy versus horizon-length trade-offs; (iii) in the case of
partial-state measurements, the issues (i) and (ii) need to be
revisited to unveil how multiple output samples from each realization
can be leveraged to recover the corresponding full state at time~$T$.

\subsection{Fixed horizon dynamic DROP}\label{subsec:prob:form}
	
We depart here from the static DROP paradigm and formulate the
\textit{dynamic} DROPs considered in this paper, where the data
evolves according to the dynamics
\begin{align}\label{data:dynamics}
  \dot{\xi}_t=F(t,\xi_t), \qquad \xi_t\in\Rat{d},
\end{align}
and is measured through the output map 
\begin{align}\label{data:output}
\zeta_t=H(t,\xi_t), \qquad \zeta_t\in\Rat{m}.
\end{align}
The initial condition $\xi_0$ is considered random with an unknown
distribution $P_{\xi_0}$. Based on the evolution of
\eqref{data:dynamics} over an horizon~$T$, we are interested in
solving a DROP with respect to the unknown distribution~$P_{\xi_T}$. As in 
Section~\ref{subsec:motivating:example}, cf.
Figure~\ref{fig:output:samples:from:trajectories}, the samples are assumed to 
be gathered from independent realizations of the system. Because of the partial 
measurements, we make the following hypothesis on how trajectories are sampled.

\begin{assumption}\longthmtitle{Sampling schedule}\label{assumption:sampling}
  The data are gathered from $\bar N$ independent trajectories
  of~\eqref{data:dynamics} over the horizon $[0,T]$, denoted by
  $\xi^i$, $i=1,\ldots,\bar N$, having i.i.d. initial conditions
  $\xi_0^i$. From each trajectory $\xi^i$ we collect $\ell$ output
  samples $\zeta_{t_i^l}^i= H(t_i^l,\xi_{t_i^l}^i)$, at $t_i^l$,
  $l\in[1:\ell]$, with $0\le t_i^1<\ldots<t_i^{\ell}\le T$ and assume
  that $0\le t_1^{\ell}\le \cdots\le t_{\bar N}^{\ell}=T$. We also
  denote $t_i^{\ell}\equiv t_i$, and, if $H(t,\xi)\equiv\xi$, we take
  $\ell=1$.
\end{assumption}

\begin{figure}[tbh]
\centering
\includegraphics[width=.9\columnwidth]{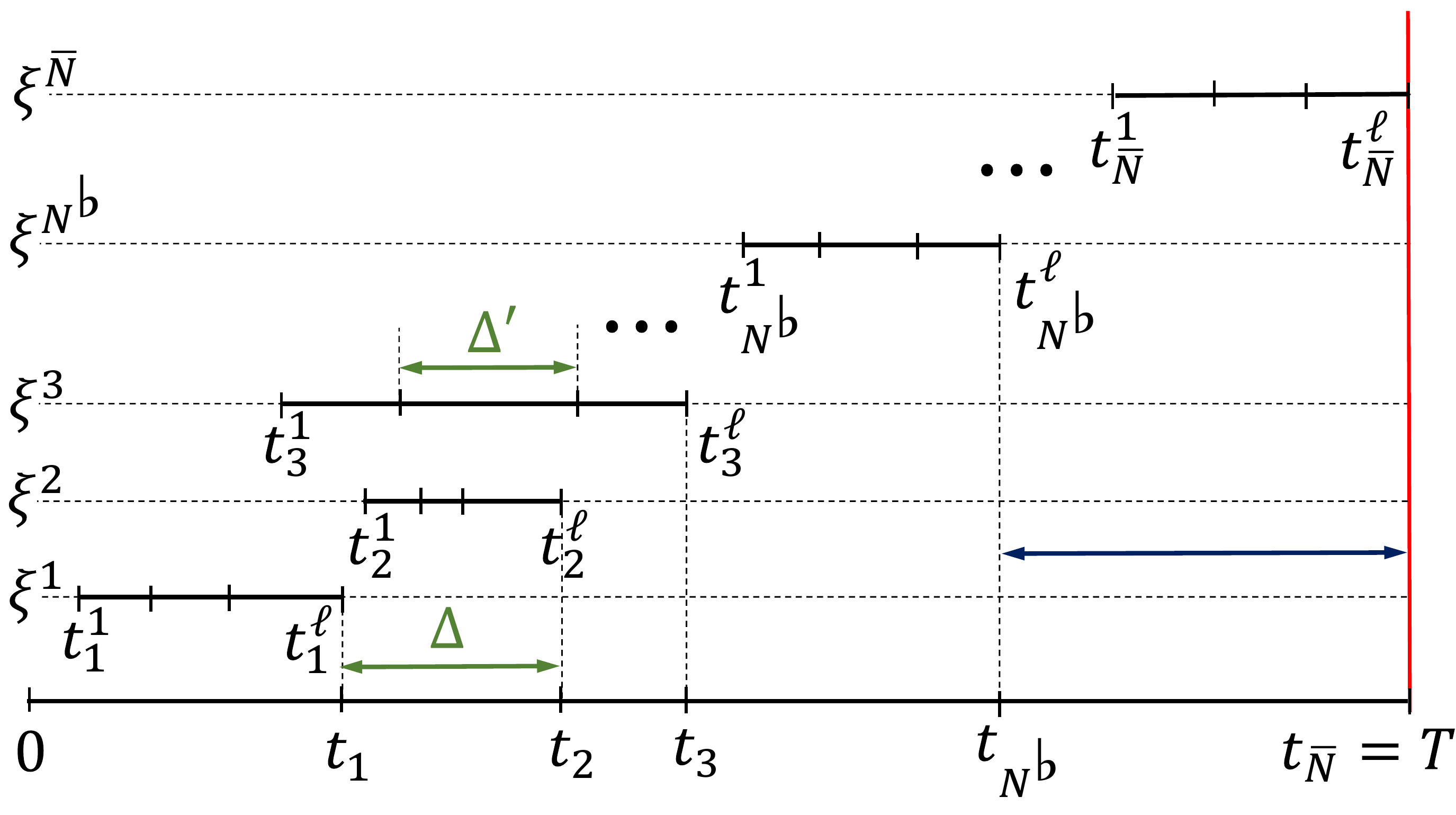}	
\caption{Sampling times where data are collected from independent trajectories 
according to Assumption~\ref{assumption:sampling}. The green arrows depict the 
tightest inter-trajectory sampling-time bound $\Delta$ and intra-trajectory 
inter-sampling-time bound $\Delta'$. The blue double-sided arrow illustrates 
the time span associated to the effective sampling horizon.}
\label{fig:sampling:seq}
\end{figure}

The assumption is motivated by scenarios where the trajectories'
outputs can be measured for a restricted amount of time. That can
happen for instance when measurements are taken by a limited-range
sensor, due to the fact that the location of the sensor and/or the
random variable's realizations changes with time
(cf. Section~\ref{sec:example}). Figure~\ref{fig:sampling:seq}
illustrates the sampling times associated to
Assumption~\ref{assumption:sampling}. We refer to $\ell$, which
represents a sufficient number of output samples to reconstruct the
full state of a trajectory, as the \textit{observation horizon
  length}.  The hypothesis that $\ell$ is common for all trajectories
is made without loss of generality, since we can otherwise take into
account only the ones from which at least $\ell$ samples are
collected, and select the last $\ell$ of them. The same holds also for
the assumption that $0\le t_1^{\ell}\le \cdots\le t_{\bar
  N}^{\ell}=T$, which can be enforced by simply relabeling the
trajectories.
 
%
%

%
We call any $\Delta>0$ such that $\Delta
\ge\max\{t_i-t_{i-1}\,|\,i\in[2:\bar N]\}$, an
\textit{inter-trajectory sampling-time bound}, and any $\Delta'>0$
such that
$\Delta'\ge\max\{t_i^l-t_i^{l-1}\,|\,l\in[2:\ell],i\in[1:\bar N]\}$ 
an \textit{intra-trajectory inter-sampling-time bound}.  In addition, we
introduce the \textit{effective sampling horizon} $[N^{\flat}:\bar N]$,
%
%
to indicate that the data used for the ambiguity set construction is
collected from the trajectories $\xi^{N^{\flat}}$ up to $\xi^{\bar N}$,
%
%
which we call \textit{effective trajectories} (cf. 
Figure~\ref{fig:sampling:seq}). The samples are exploited to estimate the state 
values of the $N:=\bar N-N^{\flat}+1$ effective trajectories at the horizon
end~$T$. The reason for considering this subset of trajectories is that 
it may not be desirable to use measurements from trajectories where
the assimilation starts before $t_{N^{\flat}}$, due to the errors
induced by the pushforward, which accumulate over time.

\textit{Problem statement:} Given the horizon $[0,T]$, and under
Assumption~\ref{assumption:sampling}, we seek to use each output tuple
$(\zeta_{t_i^1}^i,\ldots,\zeta_{t_i^{\ell}}^i)$ to estimate the state
$\widehat{\xi}_{t_i^{\ell}}^i$
of each effective trajectory $\xi^i$ at $t_i^{\ell}(\equiv t_i)$, and
determine an ambiguity set $\widehat{\P}_T^N$ containing the true
distribution $P_{\xi_T}$ with high confidence. Under a fixed
inter-trajectory sampling-time bound, we also seek to characterize the
effect of the horizon length $T$ on the size of the constructed
ambiguity set. Finally, in the presence of numerical errors or
disturbances in the dynamics, we aim to quantify the effective horizon
length up to which the ambiguity set is guaranteed to improve with the
number of samples.

We start by addressing these problems for the case of full-state
measurements, characterizing the properties of dynamic ambiguity sets
first under perfect knowledge of the flow map in
Section~\ref{sec:dyn:ambiguity} and then studying data-assimilation
versus precision trade-offs in the presence of computational errors in
Section~\ref{sec:errors:disturbances}.
Section~\ref{sec:partial:observations} extends the results to the case
when the random variable is only partially measured and evolves
linearly.

\section{Ambiguity Sets Under Full-State
  Measurements}\label{sec:dyn:ambiguity}

In this section, we characterize ambiguity sets that appear in dynamic
DROPs, as described in the problem formulation of
Section~\ref{subsec:prob:form}, when full-state measurements
$H(t,\xi)\equiv\xi$ are available. Thus, according to
Assumption~\ref{assumption:sampling}, there is an
increasing sequence of sampling times $0\le t_1\le \cdots\le t_{\bar
  N}=T$, where at each $t_i$, the state sample $\xi_{t_i}^i$ is
collected from a trajectory $\xi^i$ of \eqref{data:dynamics}, and the
observation horizon is~$\ell=1$. The ambiguity set $\widehat{\P}_T^N$
is built based on the $N=\bar N-N^{\flat}+1$ last effective samples by
leveraging concentration of measure
inequalities. 
For a moment, assume that $N$ independent samples $\xi_T^i$ become
available at time $T$, and that a selected confidence $1 - \beta >0$
is chosen. Thus, \textit{an ambiguity ball} in $\P_p(\Rat{d})$ can be
constructed using the Wasserstein metric $W_p$ with center at the
empirical distribution $\widehat P_{\xi_T}^N$ and radius
$\eps_N(\beta)$. It can be shown, cf.~\cite[Theorem 3.5]{PME-DK:17},
that the true distribution $ P_{\xi_T}$ is in this ball 
 with probability at least $1-\beta$, namely,
$\bP(W_p(\widehat P_{\xi_T}^N,P_{\xi_T})\le\eps_N(\beta))\ge 1-\beta$.
Recalling that in our setting the samples $\xi_{t_i}^i$,
$i\in[N^{\flat}:\bar N]$ are to be collected progressively prior
to $t_{\bar N}=T$, we alternatively seek to build the
\textit{cumulative empirical distribution} $\bar P_{\xi_T}^N$, using
predicted values $\bar{\xi}_T^i$ of these samples at $T$.
	
More formally, consider a probability space $(\Omega,\F,\bP)$ and a
finite sequence of i.i.d.  $\Rat{d}$-valued random variables
$(\xi_0^i)_{i\in[1:\bar N]}$ with law $P_{\xi_0}\equiv P$, where each
$\xi_0^i$ represents the initial condition of a trajectory
$\xi^i$. The trajectories evolve according to the
dynamics~\eqref{data:dynamics}. From now on, we assume that $F$ in
\eqref{data:dynamics} is continuous and locally Lipschitz in $\xi$,
and that the system is forward complete. Then, the flow map
$\Phi:\D_{\Phi}\to \Rat{d}$, where
$\D_{\Phi}:=\{(t,s,\xi)\in\RgeO\times\RgeO\times\Rat{d}:t\ge s\}$, is
defined by $\Phi(t,s,\xi):=\xi_t(s,\xi)$, inducing a family of maps
$\Phi_{t,s}:\Rat{d}\to\Rat{d}$.
Assuming that all trajectories start at time $s=0$ and using the
notation $\Phi_t=\Phi_{t,0}$, the state of each trajectory $\xi^i$ at
time $t\ge 0$ is given by the random variable
$\xi_t^i=\Phi_t\circ\xi_0^i$, with common law
$P_{\xi_t}=\Phi_{t\#}P$. We next show that, under perfect knowledge of
the flow map, the cumulative empirical distribution at $T$ formed by
the predicted values $\bar{\xi}_T^i:=\Phi_{T,t_i}(\xi_{t_i}^i)$,
coincides with that of all the samples gathered at time $T$. The proof
of the result is given in the Appendix.
	
\begin{lemma}\longthmtitle{Ideal pushforward of sampled
    states}\label{lemma:ideal:pushforward} 
  Consider a sequence of trajectories $\xi^i$ as in
  Assumption~\ref{assumption:sampling}, and the empirical distribution
  $\widehat P_{\xi_T}^N$ formed by the $N$ state samples of the
  effective trajectories $\xi^{N^{\flat}},\ldots,\xi^{\bar N}$ at $T$,
  i.e.,
  \begin{align}\label{emp:distribution:atT} \widehat
    P_{\xi_T}^N:=\frac{1}{N}\sum_{i=N^{\flat}}^{\bar N}\delta_{\xi_T^i}.
  \end{align}
  Then, (i) all $\xi_T^i$ are i.i.d., and (ii) if we consider the
  cumulative empirical distribution
  \begin{align}\label{cumulative:emp:distribution} \bar
    P_{\xi_T}^N:=\frac{1}{N}\sum_{i=N^{\flat}}^{\bar N}\delta_{\bar{\xi}_T^i},
  \end{align}
  with $\bar{\xi}_T^i=\Phi_{T,t_i}(\xi_{t_i}^i)$, $i\in[N^{\flat}:\bar N]$, 
  it holds that $\bar P_{\xi_T}^N=\widehat P_{\xi_T}^N$.
\end{lemma}
%
%


This result solves the issue of computing the evolving center of the
ambiguity ball. We next turn our attention to determining its radius
so as to ensure that the true distribution is contained in it with
high confidence.


\subsection{Ambiguity radius for compactly supported 
distributions}\label{sub:sec:concentration}

We next present results from concentration of measure to determine the
radius of the ambiguity set $\widehat{\P}_T^N$ that contains the true
distribution $P_{\xi_T}$ of the data at $T$ with a selected
confidence. Our focus is on the class of compactly supported
distributions, which is preserved under the flow of forward complete
systems. The following result
provides a concentration inequality for such laws and its explicit
dependence on the distribution's support.
	
\begin{prop}\longthmtitle{Concentration
    inequality}\label{prop:compact:support:measure}
  Consider a sequence $(X_i)_{i\in\bN}$ of i.i.d. $\Rat{d}$-valued
  random variables with a compactly supported law $\mu$. Then, for any
  $p\ge 1$, $N\ge 1$, and $\eps>0$, it holds that
  \begin{align}
    & \bP(W_p^p(\widehat{\mu}^N,\mu)\ge \eps)\le
      \chi_N(\eps,\rho;p,d), \nonumber
    \\ 
    & \chi_N(\eps,\rho):=C
      \begin{cases}
        e^{-\frac{cN}{\rho^{2p}}\eps^2} , & {\rm if}\; p>d/2,
        \\
        e^{-cN\frac{\eps^2}{\rho^{2p}(\ln(2+\rho^p/\eps))^2}}, & {\rm
          if}\; p=d/2,
        \\
        e^{-\frac{cN}{\rho^d}\eps^{\frac{d}{p}}}, & {\rm if}\; p<d/2,
      \end{cases}\label{function:chi}
  \end{align}
  where $\widehat{\mu}^N:=\frac{1}{N}\sum_{i=1}^N\delta_{X_i}$, 
  $\rho:=\frac{1}{2}\diameter{\supp{\mu}}$, and the constants $C$ 
  and $c$ depend only on $p$ and $d$.
\end{prop}
\begin{IEEEproof}
  We employ the following fact, which can be found in
  \cite[Proposition 7.16]{CV:03} and \cite[Lemma 1]{SD-MS-RS:13}.
  
  \noindent $\triangleright$ \textit{Fact I.} Let
  $T:\Rat{d}\to\Rat{d}$ with $T(x)=\bar x+Lx$ for all $x\in\Rat{d}$,
  where $\bar x\in\Rat{d}$ and $L>0$.  Then, for any
  $\mu,\nu\in\P_p(\Rat{d})$ it holds that
  $LW_p(\mu,\nu)=W_p(T_{\#}\mu,T_{\#}\nu)$. $\triangleleft$
  
  \noindent Let $z\in\Rat{d}$ with $\|x-z\|_{\infty}\le\rho$ for all
  $x\in {\rm supp}(\mu)$ and consider the mapping
  $y=T(x):=\frac{x-z}{\rho}$ and the random variables $Y_i=T(X_i)$,
  $i\in\bN$, with law $\mu_Y=T_{\#}\mu$. Then, the $Y_i$'s are also
  i.i.d. and without loss of generality they are considered supported
  precisely on $(-1,1]^d$. We claim that
  \begin{equation}\label{Wasserstein:X:vs:Y}
    \{W_p^p(\widehat{\mu}^N_Y,\mu_Y)\ge 
    \eps\}=\{W_p^p(\widehat{\mu}^N,\mu)\ge\rho^p\eps\},
  \end{equation}
  where $\widehat{\mu}_Y^N=\sum_{i=1}^N\delta_{Y_i}$. Indeed, let
  $\omega\in\Omega$.  
  Then, we have that $Y_i(\omega)=T(X_i(\omega))$
  for all $i$ and using the properties of the pushforward map, we get
  \begin{align*}
    T_{\#}\widehat{\mu}^N(\omega) &
    =T_{\#}\frac{1}{N}\sum_{i=1}^N\delta_{X_i(\omega)}
    =\frac{1}{N}\sum_{i=1}^NT_{\#}\delta_{X_i(\omega)} \\
    & =\frac{1}{N}\sum_{i=1}^N\delta_{T(X_i(\omega))}
    =\frac{1}{N}\sum_{i=1}^N\delta_{Y_i(\omega)}=\widehat{\mu}^N_Y(\omega).
  \end{align*}
  Taking also into account that $T_{\#}\mu=\mu_Y$, and exploiting
  Fact~I with $\bar x=-\frac{z}{\rho}$ and $L=\frac{1}{\rho}$, it
  follows that
  $W_p(\widehat{\mu}^N,\mu)=\rho W_p(\widehat{\mu}^N_Y,\mu_Y)$, and we
  conclude that \eqref{Wasserstein:X:vs:Y} holds. In order to derive
  the desired inequality, we will use the following result:
  
  \noindent $\triangleright$ \cite[Proposition 10]{NF-AG:15}
  \textit{Consider a sequence $(Z_i)_{i\in\bN}$ of
    i.i.d. $\Rat{d}$-valued random variables with law $\nu$, supported
    on $(-1,1]^d$.  Then, for any $p\ge 1$, $N\ge 1$, and $\eps>0$, it
    holds that
    $\bP(W_p^p(\widehat{\nu}^N,\nu)\ge \eps) \le \chi_N(\eps,1)$,
    where $\widehat{\nu}^N:=\frac{1}{N}\sum_{i=1}^N\delta_{Z_i}$, $C$
    and $c$ depend only on $p$ and $d$, and $\chi_N$ is given in
    \eqref{function:chi}}. $\triangleleft$
		
  \noindent 
  By applying this result with $\nu=\mu_Y$
  %
  %
  and $\widehat{\nu}^N=\widehat{\mu}_Y^N$ to bound
  $\bP(W_p^p(\widehat{\mu}^N_Y,\mu_Y)\ge \eps)$ and substituting the
  right-hand side of \eqref{Wasserstein:X:vs:Y} in this probability,
  we obtain the desired result.
\end{IEEEproof}

The following corollary to
Proposition~\ref{prop:compact:support:measure} characterizes the
radius of the ambiguity balls in terms of the selected confidence and
the support of the unknown distribution. 
	
\begin{corollary}\longthmtitle{Ambiguity
    radius}\label{cor:varying:compact:ambiguity}
  Under the assumptions of
  Proposition~\ref{prop:compact:support:measure}, for any confidence
  $1-\beta$, $\beta\in(0,1)$, it holds that
  $\bP(W_p(\widehat{\mu}^N,\mu)\le\eps_N(\beta,\rho))\ge 1-\beta$,
  where
  \begin{align}\label{epsN:dfn}
    \eps_{N}(\beta,\rho):=
    \begin{cases}
      \left(\frac{\ln(C\beta^{-1})}{c}\right)^{\frac{1}{2p}}\frac{\rho}{N^{\frac{1}{2p}}},
      & {\rm if}\; p>d/2,
      \\
      h^{-1}\left(\frac{\ln(C\beta^{-1})}{cN}\right)^{\frac{1}{p}}\rho,
      &
      {\rm if}\; p=d/2,
      \\
      \left(\frac{\ln(C\beta^{-1})}{c}\right)^{\frac{1}{d}}\frac{\rho}{N^{\frac{1}{d}}},
      & {\rm if}\; p<d/2,
    \end{cases}
  \end{align}
  with $h^{-1}$ the inverse of $h(x)=\frac{x^2}{(\ln(2+1/x))^2}$, $x>0$. 
\end{corollary}
%
	
\subsection{Growth conditions for ambiguity radius
  convergence}\label{subsec:dynamics:growth}
	
Here, we present sufficient conditions on the system's dynamics and
the sampling rate to guarantee that, for any prescribed confidence
$1-\beta$, the radius of the ambiguity balls converges to zero as the
horizon $[0,T]$ grows.  We start with a Lyapunov-type characterization
of the growth rate of the system dynamics.
	
\begin{prop}\longthmtitle{Lyapunov-type growth rate condition}
  \label{prop:dynamics:growth}
  For system~\eqref{data:dynamics}, assume
  that there exist a locally integrable function
  $\alpha:\RgeO\to\Rat{}$, and a function $V\in C^1(\Rat{d};\Rat{})$,
  with 
  \begin{subequations}
    \begin{align}
      a_1\|\xi\|^r\le
      & V(\xi)\le a_2\|\xi\|^r,\quad\forall\xi\in\Rat{d}, 
      \label{sandwitch:bounds}
      \\
      DV(\xi)F(t,\xi)\le & \alpha(t)V(\xi)+M_1V(\xi)^q, \nonumber
      \\
      & \hspace{5em} \forall t\ge 0,\xi\in\Rat{d}\setminus\{0\}, 
      \label{Lyapunov:inequality} 
    \end{align}
  \end{subequations}
  for certain $a_1,a_2>0$, $r>1$, $M_1\ge 0$, and
  $q\in(-\infty,1)$. Then, for any initial condition
  $\xi_0\in\Rat{d}$:
  \noindent \textit{(i)} if $M_1=0$, then
  \begin{align}\label{decay:M1:eq:zero} 
    \|\xi(t)\|\le 
    (a_2/a_1)^{\frac{1}{r}}\|\xi_0\|e^{\frac{1}{r}\int_0^t\alpha(s)ds},\quad\forall 
    t\ge 0;
  \end{align}  
  
  \noindent \textit{(ii)} if $M_1>0$, and additionally 
  \begin{align}
    \int_{t_1}^{t_2}\alpha(t)dt\le
    & M_2,\quad\forall\; t_2\ge t_1\ge 0, 
    \label{alpha:int:hypothesis}
  \end{align}
  for certain $M_2>0$, then 
  \begin{align}\label{xi:growth:rate}
    \|\xi(t)\|\le \bar M(1+\bar ct)^{\frac{1}{r(1-q)}},\quad\forall t\ge 0,
  \end{align}  
  with 
  \begin{align}\label{constant:K}
    \bar M:=(e^{M_2}(1+a_2\|\xi_0\|^r)/a_1)^{\frac{1}{r}}, \quad\bar 
    c:=M_1(1-q).
  \end{align}
\end{prop}
The proof of Proposition~\ref{prop:dynamics:growth} is given in the
Appendix.  Next, we provide the main result of this section, which
shows that for bounded inter-sampling times and under the assumption
that full-state measurements are gathered and pushed without errors
forward in time, the ambiguity sets formed by dynamics which satisfy
the growth conditions of Proposition~\ref{prop:dynamics:growth} shrink
as the DROP's horizon increases.
	
\begin{prop}\label{prop:amb:set:convergence}\longthmtitle{Ambiguity
    radius convergence} Assume that system~\eqref{data:dynamics}
  satisfies the assumptions of Proposition~\ref{prop:dynamics:growth}
  and that $P_{\xi_0}$ is supported on the compact set $K$. Select a
  confidence $1-\beta$, an exponent $p\ge 1$, and assume that
  \begin{subequations}
    \begin{align}\label{exponent:relation}
      M_1>0,\quad r(1-q)>\max\{2p,d\},
    \end{align}
    where $r,q,M_1$ are given in
    Proposition~\ref{prop:dynamics:growth}, or that
    \begin{align}\label{exponent:relation:alternative}
      M_1=0,\quad rq'>\max\{2p,d\},\quad \int_0^t\alpha(s)ds\le \frac{\ln 
        t}{q'},\forall t\ge t_0
    \end{align}
  \end{subequations}
  for some $q'>0$ and $t_0 >0$. 
  For any horizon $[0,T]$, consider a sampling sequence as in
  Assumption~\ref{assumption:sampling}, with $t_1=0$ and
  inter-trajectory sampling-time bound independent of $T$, and let
  \begin{align}\label{rho:T}
    \rho_T:=\diameter{\Phi_T(K)}/2,
  \end{align}
  with $\Phi_T(K)$ the set of reachable states at $T$ from $K$. Then,
  \begin{subequations}
    \begin{align} 
      P(W_p(\bar P_{\xi_T}^N,P_{\xi_T})\le\eps_N({\beta},\rho_T)) & \ge 
      1-\beta, \label{ambig:set:dfn:and:convergence1} \\ 
      \lim_{T\to\infty} \eps_N({\beta},\rho_T) & 
      =0,\label{ambig:set:dfn:and:convergence2}
    \end{align}  
  \end{subequations}
  %
  %
  where $\eps_N$ is given in~\eqref{epsN:dfn} and $\bar P_{\xi_T}^N$
  is the cumulative empirical distribution in
  \eqref{cumulative:emp:distribution}, with $N^{\flat}=1$ and $N=\bar
  N$.
\end{prop}
\begin{IEEEproof}
  We will leverage the results of Lemma~\ref{lemma:ideal:pushforward},
  Corollary~\ref{cor:varying:compact:ambiguity}, and
  Proposition~\ref{prop:dynamics:growth} for the proof. Note first
  that, under perfect knowledge of the flow map, it follows from
  Lemma~\ref{lemma:ideal:pushforward}(ii) that the center $\bar
  P_{\xi_T}^N$ of the ambiguity ball in
  \eqref{ambig:set:dfn:and:convergence1} is the same as $\widehat
  P_{\xi_T}^N$, i.e., the empirical distribution formed by $N$ samples
  $\xi_T^1,\ldots,\xi_T^N$ of independent trajectories at $T$. Due to
  Lemma~\ref{lemma:ideal:pushforward}(i), these samples are
  i.i.d.. Thus, we can apply
  Corollary~\ref{cor:varying:compact:ambiguity} to the empirical
  distribution $\widehat{\mu}^N=\widehat P_{\xi_T}^N = \bar
  P_{\xi_T}^N$ to infer that $\bP(W_p(\bar
  P_{\xi_T}^N,P_{\xi_T})\le\eps_N(\beta,\rho_T))\ge 1-\beta$ holds,
  or, in other words,
  ~\eqref{ambig:set:dfn:and:convergence1}.

  Next, note that since $K$ is compact, there is some $\rho>0$ with
  $K\subset \{\xi\in\Rat{d}:\|\xi\|_{\infty}\le\rho\}$. Thus, when
  \eqref{exponent:relation} is fulfilled, it follows from
  Proposition~\ref{prop:dynamics:growth} that
  \begin{align}\label{rhoT:bund}
    \rho_T\le \bar M(1+\bar cT)^{\frac{1}{r(1-q)}},
  \end{align}
%
  %
  where $\bar c$ and $\bar M$ are given by \eqref{constant:K} with
  $\|\xi_0\|$ in the definition of $\bar M$ replaced by $\sqrt d\rho$.
  %
  %
  Next, define $i(T):={\lfloor T/\Delta \rfloor}+1$, where
  $\Delta$ is a common inter-trajectory sampling-time bound for each
  horizon $[0,T]$. Note that $i(T)\le N$ and $T<\Delta i(T)$,
  implying by \eqref{rhoT:bund} that
  \begin{align*}
    \rho_T<\bar M(1+\bar c\Delta
    i(T))^{\frac{1}{r(1-q)}}=:\bar{\rho}_T.
  \end{align*}
  From the latter, \eqref{epsN:dfn}, the fact that $i(T)\le N$, and
  that $\eps_N(\beta,\rho_T)$ decreases with $N$ and increases with
  $\rho_T$, it follows that
  \begin{align*}
    & \eps_N(\beta,\rho_T)<\eps_{i(T)}(\beta,\bar{\rho}_T) \\
    & =\begin{cases}
      \left(\frac{\ln(C\beta^{-1})}{c}\right)^{\frac{1}{2p}}\frac{\bar
        M(1+\bar c\Delta
        i(T))^{\frac{1}{r(1-q)}}}{i(T)^{\frac{1}{2p}}}, & {\rm if}\;
      p>d/2, \\
      h^{-1}\left(\frac{\ln(C\beta^{-1})}{ci(T)}\right)^{\frac{1}{p}}\bar
      M(1+\bar c\Delta i(T))^{\frac{1}{r(1-q)}}, & {\rm if}\; p=d/2, \\
      \left(\frac{\ln(C\beta^{-1})}{c}\right)^{\frac{1}{d}}\frac{\bar
        M(1+\bar c\Delta i(T))^{\frac{1}{r(1-q)}}}{i(T)^{\frac{1}{d}}}
      & {\rm if}\; p<d/2.
    \end{cases}
  \end{align*}
  Thus, it suffices to show that
  $\lim_{T\to\infty}\eps_{i(T)}(\beta,\bar{\rho}_T)=0$. In particular,
  when $p \neq d/2$, by setting $\bar p=\max\{2p,d\}$, we get
  \begin{align*}
    \lim_{T\to\infty}\eps_{i(T)}(\beta,\bar{\rho}_T) &
    =\lim_{T\to\infty}\frac{C'(1+\bar c\Delta
      i(T))^{\frac{1}{r(1-q)}}}{i(T)^{\frac{1}{\bar p}}} \\
    & \le \lim_{T\to\infty}C''i(T)^{\frac{\bar p-r(1-q)}{\bar
        pr(1-q)}}=0,
  \end{align*}
  because of \eqref{exponent:relation} and the fact that
  $i(T)\to\infty$ when $T\to\infty$, where the constants $C'$ and
  $C''$ in the derivation are independent of $T$. For the case where
  $p=d/2$, the result follows analogously by exploiting the following
  technical fact whose proof is given in the Appendix.

  \noindent $\triangleright$ \textit{Fact II.} For any $\bar q>2$ and
  $a>0$, it holds that
  $\lim_{\kappa\to\infty}h^{-1}(a/\kappa)\kappa^{\frac{1}{\bar q}}=0$,
  with $h$ given in Corollary~\ref{cor:varying:compact:ambiguity}.
  $\triangleleft$

  Finally, when \eqref{exponent:relation:alternative} holds, we obtain
  form \eqref{decay:M1:eq:zero} that
  \begin{align*}
    \|\xi(t)\|\le (a_2/a_1)^{\frac{1}{r}}\|\xi_0\|e^{\frac{1}{rq'}\ln
      t},
  \end{align*}
  for all $t\ge t_0$, and thus, that there exists a constant $\bar M'$
  with $\rho_T\le \bar M'(1+T)^{\frac{1}{rq'}}$ for all $T>0$. Then,
  the result follows by using precisely the same arguments as before.
\end{IEEEproof}	     
	

\section{Ambiguity Sets under Pushforward Errors and
  Disturbances}\label{sec:errors:disturbances}
	
In general, the dynamics~\eqref{data:dynamics} of the random variable
cannot be solved for explicitly. The pushforward of full-state samples
hence requires the numerical integration of the system's solutions,
which gives rise to an approximation~$\supscr{\Phi}{num}$ of the exact
flow map.
In addition, the dynamics of the process may be subject to
disturbances. Both cases suggest that the cumulative empirical
distribution $\bar P_{\xi_T}^N$ in~\eqref{cumulative:emp:distribution}
will no longer coincide with the empirical
distribution~$\widehat P_{\xi_T}^N$
in~\eqref{emp:distribution:atT}. We proceed to quantify this
difference and characterize the relevant ambiguity sets, focusing
first on the case where the flow is numerically integrated.
%
%
%
Recall that the dynamics is continuous and locally Lipschitz. It 
therefore follows from classical approaches to bound the 
numerical integration error~\cite[Theorem 3.4.7]{AS-ARH:98}, that for any 
compact set $K\subset\Rat{n}$ and finite
time horizon $[0,T]$, there exist positive constants $\fK$ and $L$,
such that
\begin{align}\label{numerical:error:bound}
\|\supscr{\Phi}{num}_{t,s}(\xi)-\Phi_{t,s}(\xi)\|\le
\fK(e^{L(t-s)}-1),
\end{align}
for all $0\le s\le t\le T$ and $\xi\in\Phi_s(K)$. For the center of
the ambiguity ball containing $P_{\xi_T}$, we consider a variant of
the cumulative empirical distribution $\bar P_{\xi_T}^N$ in
\eqref{cumulative:emp:distribution}, formed by pushing forward the
full-state samples by~$\supscr{\Phi}{num}$.  The following result,
whose proof is in the Appendix, specifies the radius of the ambiguity
set at the end of the horizon.

\begin{thm}\longthmtitle{Ambiguity
    radius with approximate pushforward}\label{thm:numerical:ambiguity:radius}
  Assume that the support of the initial condition of
  system~\eqref{data:dynamics} is contained in the compact set~$K$.
  Consider a sampling sequence as in
  Assumption~\ref{assumption:sampling} with inter-trajectory
  sampling-time bound $\Delta>0$, and the cumulative empirical
  distribution $\bar P_{\xi_T}^N$ in
  \eqref{cumulative:emp:distribution}, with
  $\bar{\xi}_T^i:=\supscr{\Phi}{num}_{T,t_i}(\xi_{t_i}^i)$,
  $i\in[N^{\flat}:\bar N]$. Then, for any confidence $1-\beta$ and
  $p\ge 1$ it holds
  \begin{subequations}
    \begin{align}
      \bP(W_p(\bar P_{\xi_T}^N,P_{\xi_T})& \le\psi_N(\beta,\Delta))\ge
      1-\beta,
      \label{total:ambituity:characterization}
      \end{align}
      where
      \begin{align}
      \psi_N(\beta,\Delta): &
      =\eps_N(\beta,\rho_T)+\bar{\eps}_N(\Delta),
      \label{mixed:ambiguity:radius} 
      \\
      \bar{\eps}_N(\Delta): & 
      = \fK \Big(\frac{1}{N}\int_1^N(e^{L\Delta s}-1)^pds
      \Big)^{\frac{1}{p}},
      \label{bar:epsN}
    \end{align}
  \end{subequations}
  with $\eps_N(\beta,\rho_T)$ and $\rho_T$, as given in
  \eqref{epsN:dfn} and \eqref{rho:T}, respectively.
\end{thm}


According to this result, when the exact flow map is no longer
available, the ambiguity radius needs to be increased with the
additional term $\bar \eps_N$ to obtain the same high-confidence
guarantee. In the ideal case where $\fK=0$ in~\eqref{bar:epsN}, we
recover the result of
Corollary~\ref{cor:varying:compact:ambiguity}. The extra term
$\bar \eps_N$ increases as we consider more samples, since they are
located farther back in time and pushing them forward induces larger
errors.  Based on Theorem~\ref{thm:numerical:ambiguity:radius} we
quantify the effective sampling horizon size in terms of guaranteed
ambiguity reduction in the presence of numerical errors.
	
\begin{prop}\longthmtitle{Effective sampling
    horizon}\label{prop:effective:horizon}
  Under the hypotheses of
  Theorem~\ref{thm:numerical:ambiguity:radius}, assume
  additionally that $p\ne d/2$. Let $\bar p:=\max\{2p,d\}$ and
  $\bar C\equiv\bar C(\bar
  p,\beta,\rho_T):=\Big(\frac{\ln(C\beta^{-1})}{c}\Big)^{\frac{1}{\bar
      p}}\rho_T$, with $C$, $c$ and $\rho_T$ as in \eqref{epsN:dfn}
  and \eqref{rho:T}, respectively. Then, there exists
  $\Delta^*\equiv\Delta^*(\bar C,L,\fK)>0$ such that for every
  $\Delta\in(0,\Delta^*)$, the set
  \begin{align*}
    \N(\Delta):=\{N\in\bN\,|\,\bar C(
    & \kappa^{-\frac{1}{\bar       
      p}}-(\kappa+1)^{-\frac{1}{\bar p}}) \\
    & >\bar{\eps}_{\kappa+1}(\Delta)-\bar{\eps}_{\kappa}(\Delta),\,
      \forall\kappa\in[1:N]\},
  \end{align*}
  with $\bar{\eps}_{1},\ldots, \bar{\eps}_{N+1}$ as given in
  \eqref{bar:epsN},
  is nonempty and bounded. In addition, the ambiguity radius
  $\psi_N(\beta,\Delta)$ in \eqref{mixed:ambiguity:radius} is strictly
  decreasing with $N$, for $N\in[1:N^*(\Delta)]$, with
  $N^*(\Delta):=\max(\N(\Delta))+1$.
\end{prop}
\begin{IEEEproof}
  From \eqref{bar:epsN}, we get that
  $\lim_{\Delta\to 0}\bar{\eps}_N(\Delta)=0$ for all $N\in\bN$. Thus,
  there exists $\Delta^*>0$ such that
  $\bar{\eps}_{2}(\Delta)-\bar{\eps}_1(\Delta)<\bar
  C(1-2^{-\frac{1}{\bar p}})$ for all $\Delta<\Delta^*$, implying that
  $1\in\N(\Delta)$, and hence, that $\N(\Delta)$ is nonempty. In
  addition, we have that
  $\lim_{N\to\infty} \bar{\eps}_N(\Delta)=\infty$ for all $p\ge 1$,
  which implies that $\N(\Delta)$ is bounded. Indeed, otherwise it
  would hold that
  \begin{align*}
    & \bar{\eps}_{N+1}(\Delta)-\bar{\eps}_1(\Delta)  
      =\sum_{\kappa=1}^N(\bar{\eps}_{\kappa+1}(\Delta) 
      -\bar{\eps}_{\kappa}(\Delta)) \\
    & <\sum_{\kappa=1}^N\bar C(\kappa^{-\frac{1}{\bar 
      p}}-(\kappa+1)^{-\frac{1}{\bar p}}) 
      =\bar C(1-(N+1)^{-\frac{1}{\bar p}})<\bar C 
  \end{align*}
  for all $N\in\bN$, leading to a contradiction.  Thus $N^*(\Delta)$
  is always finite. Furthermore, from \eqref{epsN:dfn}, 
  \eqref{mixed:ambiguity:radius}, and the definition
  of $\bar C$, $N^*(\Delta)$, and $\N(\Delta)$, we get 
  $\psi_{N+1}(\beta,\Delta)<\psi_N(\beta,\Delta)$ for all
  $N\in[1: N^*(\Delta) -1]$, as desired.
\end{IEEEproof}

\begin{figure}[tbh]
  \centering
  \subfigure[]{\includegraphics[width=.55\linewidth]{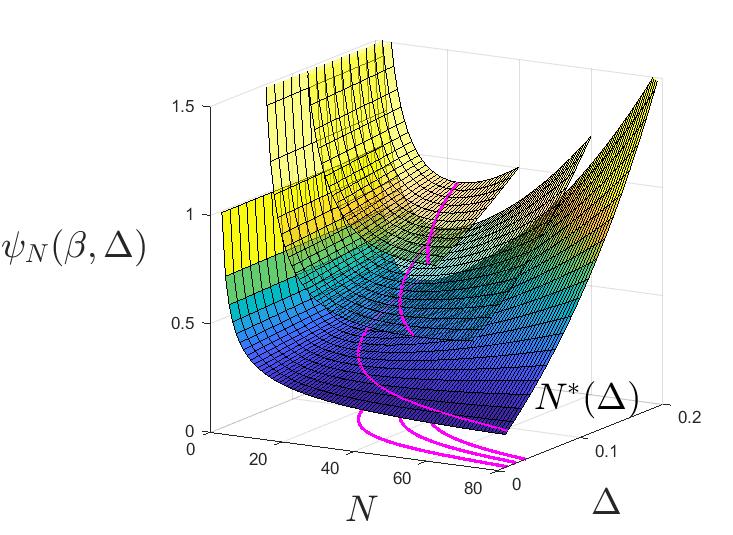}}
  \hspace*{-2ex}
  \subfigure[]{\includegraphics[width=.45\linewidth]{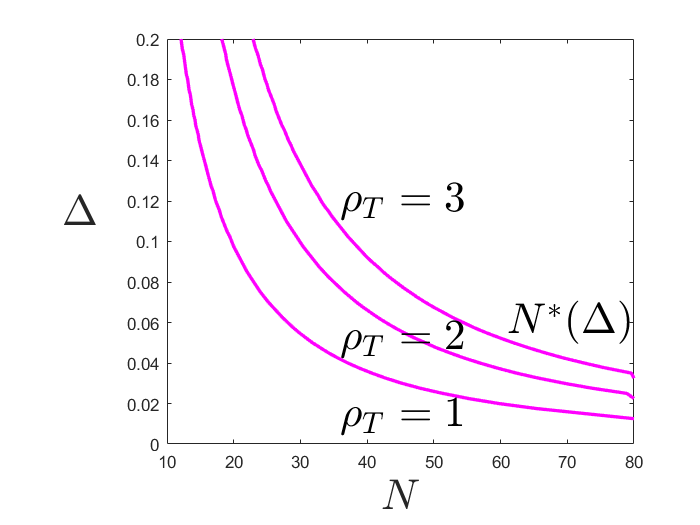}}
  \caption{(a) shows how the ambiguity radius $\psi_N(\beta,\Delta)$
    in~\eqref{mixed:ambiguity:radius} varies with respect to $\Delta$
    and $N$ {for a fixed confidence $1-\beta$, and the parameter
      values $d=1$, $p=1$, $L=0.1$, $\fK=1$,} and $\rho_T \in
    \{1,2,3\}$.  Given that the component $\eps_N$ of the radius,
    which is strictly decreasing with $N$, is proportional to the
    distributions' support size and that the effect of numerical
    errors is independent of $\rho_T$, the effective sampling horizon
    increases with $\rho_T$, as shown in (b).}
  \label{fig:amiguity:radius}
\end{figure}

Proposition~\ref{prop:effective:horizon} identifies the upper bound
$N^*(\Delta)$ on the number of samples below which one can guarantee
that the more samples the better regarding the ambiguity radius.
Figure~\ref{fig:amiguity:radius} illustrates how the ambiguity radius
$\psi_N(\beta,\Delta)$ in~\eqref{mixed:ambiguity:radius} varies with
respect to the sampling period and the number of samples in the
presence of numerical errors.

%

Theorem~\ref{thm:numerical:ambiguity:radius}
and Proposition~\ref{prop:effective:horizon} are also applicable when 
the sampled trajectories are subject to unknown disturbances~$\mathfrak
d$. Formally, the dynamics in this case takes the form
\begin{align*}
  \dot{\xi_t}=F(t,\xi_t,\mathfrak d_t) ,
\end{align*}
where $\mathfrak d$ belongs to a class $\mathcal D$ of functions, which
are uniformly bounded for every finite time horizon, and with $\mathfrak 
d\equiv 0$ being an element of $\mathcal D$. Additionally, we
assume that $F$ is also locally Lipschitz with respect to its
$\mathfrak d$ argument.  For any constant $\epsilon>0$ and horizon
$[0,T]$, let $B:=\{\Phi_t(\xi)\,|\,\xi\in K,t\in[0,T]\}$ and
$B_\epsilon:=\{\xi\in\Rat{d}\,|\,\dist(\xi,B)\le \epsilon\}$, where
$K$ is a compact set containing the initial state~$\xi_0$.  Using the
local Lipschitzness assumption, we can select a constant $L>0$ so that
\begin{align*}
  \|F(t,\xi,\mathfrak d_t)-F(t,\xi',\mathfrak d_t')\|\le L\|(\xi,\mathfrak 
  d_t)-(\xi',\mathfrak d_t')\|,
\end{align*} 
for all $t\in [0,T]$, $\xi,\xi'\in B_\epsilon$ and
$\mathfrak d,\mathfrak d'\in \mathcal D$. Then, from~\cite[Theorem 
3.4]{HKK:02},
it follows that~\eqref{numerical:error:bound} holds with
$\supscr{\Phi}{num}_{t,s}$ replaced by the flow of
$\dot{\xi_t}=F(t,\xi_t,\mathfrak d_t)$ and
$\fK:=\min\{\sup_{t\in[0,T],\mathfrak d\in\mathcal D}\mathfrak
d_t,\frac{\epsilon}{e^{LT}-1}\}$, for all disturbances $\mathfrak 
d\in\mathcal D$ with $\max_{t\in[0,T]}\mathfrak d_t\le\fK$. With this bound 
in place, the proofs of Theorem~\ref{thm:numerical:ambiguity:radius} and 
Proposition~\ref{prop:effective:horizon} also hold for the case with
disturbances.
	
\section{Ambiguity Sets for Partially Observable Linear
  Systems}\label{sec:partial:observations}
	
Here we consider the case of partial measurements, where multiple
samples are collected from each independent trajectory.  We restrict
our attention to linear time-varying systems with linear outputs,
i.e.,
\begin{align}\label{linear:system}
  \dot \xi_t=A(t)\xi_t,\quad \zeta_t=C(t)\xi_t,
\end{align}
with $A(t)\in\Rat{d\times d}$ and $C(t)\in\Rat{m\times d}$. In this
case, the full state of each trajectory $\xi^i$ is no longer directly
available through the individual output samples. However, it is
possible to recover it when the system is sampled-data observable, and
combine this knowledge with our approach in the previous sections to
build the ambiguity sets. We therefore start by examining
observability conditions under which the state reconstruction is
possible when a sufficient number of output samples is collected.

\subsection{Sampled data observability}
For the linear time-varying system \eqref{linear:system}, let
$\Phi(t,s)\in\Rat{d\times d}$, for $t,s\ge 0$, denote its fundamental
matrix, satisfying $\Phi_{t,s}(\xi)=\Phi(t,s)\xi$ for all
$\xi\in\Rat{d}$.  According to Assumption~\ref{assumption:sampling},
we have $\ell$ output samples from each trajectory $\xi^i$.  Each such
sample can be evaluated as
$\zeta_{t_i^l}^i=C(t_i^l)\Phi(t_i^l,t_i^{\ell})\xi_{t_i^\ell}^i$,
$l\in[1:\ell]$, by taking the state at $t_i^{\ell}$ backward to time
$t_i^l$ through the flow and computing its output value. Consequently,
recovering the unmeasured state $\xi_{t_i^\ell}^i$ at the last
sampling instant is equivalent to requiring that the
\textit{sample-observability matrix}
\begin{align}\label{general:observability:matrix:timevar}
  \Obs_i\equiv \Obs_{t_i^1\cdots t_i^{\ell}}^{\exp-} 
  :=\left(\begin{matrix}
      C(t_i^1)\Phi(t_i^1,t_i^{\ell}) \\
      C(t_i^2)\Phi(t_i^2,t_i^{\ell})  \\
      \vdots \\
      C(t_i^{\ell}) 
    \end{matrix}\right), 
\end{align}
is left invertible. This turns out to hold when
system~\eqref{linear:system} is sampled data observable and the
observation horizon length~$\ell$ is sufficiently large. We next
present such sampled-data observability results, starting with the
case when the system is (linear) time-invariant (LTI), i.e.,
$A(t)\equiv A$ and $C(t)\equiv C$.

\begin{assumption}\label{sampled:observability:assumptions}
  The pair $(A,C)$ is observable and the sampling schedule satisfies
  either one of the following hypotheses:

  \noindent \textbf{H1} \longthmtitle{Equidistant sampling} The
  sampling times are given by $t_i^l=t_i^1+(l-1)\Delta'$,
  $l\in[1:\ell]$, $i\in[1:\bar N]$, for some $\Delta'>0$, with
  $\Delta'(\lambda-\lambda')\ne 2k\pi j$ for all $k\in\mathbb Z$ and
  distinct eigenvalues $\lambda$, $\lambda'$ of $A$, and $\ell\ge d$
  (denoting $j\equiv\sqrt{-1}\in\mathbb C$).

  \noindent \textbf{H2} \longthmtitle{Periodic non-equidistant
    sampling} The sampling times are given by the pattern
  \begin{align*}
    t_i^l:=\begin{cases}
      t_i^1+(l-1)\Delta',  & l\in[1:\bar d+1], \\
      t_i^1+\bar d\Delta'+\Delta'',  & l=\bar d+2, \\
      t_i^{l-(\bar d+1)}+\bar d\Delta'+\Delta'', & l\in[\bar
      d+3:\ell],
    \end{cases}
  \end{align*}
  with $\Delta',\Delta''>0$, $\Delta'/\Delta''\notin\mathbb Q$, $\bar
  d\ge d$, and $\ell\ge (\bar d+1)d$.

  \noindent \textbf{H3} \longthmtitle{Irregular sampling} The
  observation horizon length $\ell$ satisfies the lower bound
  \begin{align*}
    \ell>\fm-1+\tau\delta/(2\pi),
  \end{align*}
  where $\delta:=\max_{1\le j,j'\le q}\{\Im(\lambda_j'-\lambda_j)\}$,
  $\fm:=\sum_{j=1}^q\fm_j$,
  $\tau:=\max\{t_i^{\ell}-t_i^1, {i\in[1:\bar N]}\}$,
  and $\fm_1,\ldots,\fm_q$ are the indices of~$A$'s eigenvalues
  $\lambda_1,\ldots,\lambda_q$.
\end{assumption}

The first sampling hypothesis in 
Assumption~\ref{sampled:observability:assumptions} is a classical result, 
commonly known as the Kalman-Ho-Narendra criterion (see e.g.,~\cite{EDS:98}).
The second hypothesis is given in~\cite{GK:99}, 
whereas the last observability criterion under irregular sampling was
originally presented in~\cite{LYW-CL-GY-LG-CZX:11}. Here, we provide a
slightly refined version of the latter from the more recent
work~\cite[Theorem 2]{SZ-FA:16}. We next state formally that under
either of the three cases for the sampling schedule, the sample
observability matrix $\Obs_i$ is left invertible. The proof of this
result invokes classical arguments from LTI system theory. For
completeness, since we found no verbatim statement in the literature,
we provide a brief proof in the Appendix.

\begin{lemma}\longthmtitle{Sampled-data observability for LTI
    systems}\label{lemma:sampled:observability:LTI} 
  Assume that system~\eqref{linear:system} is LTI and that samples are
  collected according to Assumption~\ref{assumption:sampling}. Then,
  under the observability
  Assumption~\ref{sampled:observability:assumptions}, each matrix
  $\Obs_i$ is left invertible.
\end{lemma}

Left invertibility of $\Obs_i$ is equivalent to the property that its
smallest singular value is strictly positive. However, the result of
Lemma~\ref{lemma:sampled:observability:LTI} does not provide a bound
on how far the smallest singular value lies from zero, or in other
words, on the distance of the matrix from becoming singular. This
distance is of particular interest when the collected output samples
are perturbed.

We next provide conditions which guarantee robust invertibility of the
sample-observability matrix by deriving uniform lower bounds for the
smallest singular value of a weighted variant of~$\Obs_i$.  These
constitute a robust observability criterion under irregular sampling
for the general time-varying case~\eqref{linear:system}.
We make the following assumption:

\begin{assumption}\label{observability:assumption:tv}
  \longthmtitle{Time-varying observability} System
  \eqref{linear:system} is observable on any interval $[a,b]\subset
  [0,T]$, $t\mapsto A(t)$ is continuous, and $t\mapsto C(t)$ is
  continuously differentiable.
\end{assumption}

For each trajectory $\xi^i$, let $\tau(i):=t_i^{\ell}-t_i^1$ denote
the length of the time interval during which output samples are
collected.  It holds that
\begin{align*} 
  0<\supscr{\tau}{low} :=\min_{i\in[1:\bar N]}\tau(i) \le \tau(i) \le
  \supscr{\tau}{up} :=\max_{i\in[1:\bar N]}\tau(i) \le T,
\end{align*}
for all $i\in[1:\bar N]$.  Next, define
\begin{align}\label{map:K}
  K(s,t):=\Phi(s,t)^{\top}C(s)^{\top}C(s)\Phi(s,t)
\end{align} 
and 
\begin{align}\label{Wtau} 
  \W_{\varsigma}(t):=\int_{t}^{t+\varsigma}K(s,t+\varsigma)ds,
\end{align}
for $s,t,\varsigma\ge 0$.  Then, 
$\W_{\varsigma}(t)=\Phi({t,t+\varsigma})^{\top}
\widetilde{\W}_{\varsigma}(t)\Phi({t,t+\varsigma})$, where
$\widetilde{\W}_{\varsigma}(t):=\int_{t}^{t+\varsigma}K(s,t)ds$ is the
observability Gramian of \eqref{linear:system} on $[t,t+\varsigma]$,
which is continuous with respect to $t$. By
Assumption~\ref{observability:assumption:tv}, the observability
Gramian $\widetilde{\W}_{\supscr{\tau}{low}}(t)$ is also positive
definite for all $t\in[0,T-\supscr{\tau}{low}]$ (see e.g.,
\cite[Exercise 6.3.2]{EDS:98}). Thus, the same properties hold for
$\W_{\supscr{\tau}{low}}(t)$, implying
\begin{align}\label{minimum:eigenvalue:over:Gramians}
  \lambda_{\min}\big(\W_{\supscr{\tau}{low}}
  (t)|_0^{T-\supscr{\tau}{low}}\big)
  :=\min_{t\in[0,T-\supscr{\tau}{low}]}
  \lambda_{\min}(\W_{\supscr{\tau}{low}}(t))>0,
\end{align}
where $\lambda_{\min}$ denotes smallest eigenvalue. Next, fix any
$i\in[1:\bar N]$ and let
\begin{align}
  {\tau_l(i)}:= t_i^{l+1}-t_i^l, \;
  l\in[1:\ell-1],\label{times:hatt:tau}
\end{align}
be the lengths of the inter-sampling time intervals.
Define the weight matrix
\begin{align}\label{weight:matrix:W}
  \bm W_i & := \diag{w_1(i),\ldots,w_{\ell}(i)} \otimes I_m,
  \\
  & \quad w_1^2(i) =\frac{\tau_1(i)}{2}, \;
  w_l^2(i)=\frac{\tau_{l-1}(i)+\tau_l(i)}{2},
  \notag  
  \\
  & \hspace{3.4em} l\in[2:\ell-1],\; w_{\ell}^2(i)
  =\frac{\tau_{\ell-1}(i)}{2}.  \notag
\end{align}
%

Next, we find a positive lower bound for the smallest eigenvalue of
$\Obs_i^{\top}\bm W_i\bm W_i \Obs_i$. This fact implies the left invertibility
of~$\Obs_i$.  The proof of the result is given in the Appendix.
      
\begin{prop}\longthmtitle{Robust sampled-data
    observability}\label{prop:robust:observability}
  Under Assumption~\ref{observability:assumption:tv}, for $a\in(0,1)$,
  assume the intra-trajectory inter-sampling-time bound $\Delta'$
  satisfies
  \begin{align}\label{Delta:for:sampled:observability}
    \Delta' \le
    \frac{4(1-a)\lambda_{\min}\big(\W_{\supscr{\tau}{low}}(t)|_0^{T-\supscr{\tau}{low}}\big)}
     {\supscr{\tau}{up}\max_{\supscr{\tau}{low}\le t\le
        T,\max\{0,t-\supscr{\tau}{up}\}\le s\le t}\|K_s(s,t)\|},
  \end{align}
  where $K_s(s,t):=\frac{\partial}{\partial s}K(s,t)$. Then, the system is 
  sampled-data observable, i.e., for any $i\in[1:\bar N]$ the matrix $\Obs_i$ 
  is invertible. In addition, 
  \begin{align}\label{eigenvalue:inequality}
    \lambda_{\min}(\Obs^{\top}_i\bm W_i\bm W_i\Obs_i) \ge 
    a\lambda_{\min}\big(\W_{\supscr{\tau}{low}}(t)|_0^{T-\supscr{\tau}{low}}\big).
  \end{align} 
\end{prop}

We employ the bound in Proposition~\ref{prop:robust:observability} to
quantify the state reconstruction error under perturbed measurements
in the next section.
The result of Proposition~\ref{prop:robust:observability} takes a more
explicit form for LTI systems. In this case, note that
$K(s,t)=\widehat K(s-t)$, with
\begin{align}\label{function:hatK}
  \widehat K(t):=e^{A^{\top}t}C^{\top}Ce^{At}.
\end{align}
Hence,
$\widehat{\W}_{\varsigma}:=\W_{\varsigma}(t)=\int_0^{\varsigma}\widehat
K(s-{\varsigma})ds$ is independent of $t$, and
$\lambda_{\min}\big(\W_{\supscr{\tau}{low}}(t)|_0^{T-\supscr{\tau}{low}}\big)
= \lambda_{\min}(\widehat{\W}_{\supscr{\tau}{low}})$.

\begin{corollary}\label{corollary:robust:observability:LTI}
  \longthmtitle{Robust sampled-data observability for LTI systems}
  Under the assumptions of
  Proposition~\ref{prop:robust:observability}, when
  \eqref{linear:system} is time-invariant and the intra-trajectory
  inter-sampling-time bound $\Delta'$ satisfies
  \begin{align}\label{Delta:for:LTI:sampled:observability}
    \Delta'\le
    \frac{2(1-a) \lambda_{\min}(\widehat 
    \W_{\supscr{\tau}{low}})}{\supscr{\tau}{up}
      \max_{s\in[0,\supscr{\tau}{up}]}\|\widehat
      K(s-\supscr{\tau}{up})A\|},
  \end{align}
  %
  the system is sampled-data observable and \eqref{eigenvalue:inequality} holds.
\end{corollary} 

\subsection{Ambiguity sets under exact and inexact observations}

Here we exploit the sampled observability results to characterize
dynamic ambiguity sets under progressively collected output
measurements. The next result is the analogue of
Lemma~\ref{lemma:ideal:pushforward} when~\eqref{linear:system} is
sampled data observable and measurements are exact.

\begin{lemma}\longthmtitle{Ideal pushforward of output-sample
    reconstructed states}\label{lemma:ideal:output:pushforward}
  Consider a sequence of trajectories $\xi^i$ as in
  Assumption~\ref{assumption:sampling}, the empirical distribution
  $\widehat P_{\xi_T}^N$ in~\eqref{emp:distribution:atT}, and the
  cumulative empirical distribution $\bar P_{\xi_T}^N$ given
  by~\eqref{cumulative:emp:distribution}, with
  \begin{align}
    \bar{\xi}_T^i:= {\Phi(T,t_i^{\ell})}\Obs_i^{\dagger}
    \bm{\zeta}^i,
    \label{xiibar:atT:output:dfn}
  \end{align}
  for $i\in[N^{\flat}:\bar N]$, where $ \bm{\zeta}^i:= (
  \zeta_{t_i^1}^i , \dots, \zeta_{t_i^{\ell}}^i )$. 
  %
  %
  Assume that either~\eqref{linear:system} is time-invariant and
  Assumption~\ref{sampled:observability:assumptions} or the hypotheses
  of Corollary~\ref{corollary:robust:observability:LTI} hold, or that
  it is time-varying and the assumptions of
  Proposition~\ref{prop:robust:observability} hold. Then, $\bar
  P_{\xi_T}^N=\widehat P_{\xi_T}^N$.
\end{lemma}
\begin{IEEEproof}
  From the definition of $\widehat P_{\xi_T}^N$ and $\bar
  P_{\xi_T}^N$, the result follows if we establish that
  $\bar{\xi}_T^i=\xi_T^i$ for all $i\in[N^{\flat}:\bar N]$. For each
  of the possible cases, we get from either
  Lemma~\ref{lemma:sampled:observability:LTI},
  Proposition~\ref{prop:robust:observability}, or
  Corollary~\ref{corollary:robust:observability:LTI}, that the
  matrices $\Obs_i$ are left invertible, and hence, that
  $\Obs_i^{\dagger} \bm{\zeta}^i=\xi_{t_\ell^i}^i$. Thus, we deduce
  that $\bar{\xi}_T^i :=
  {\Phi(T,t_i^{\ell})}\xi_{t_\ell^i}^i=\xi_T^i$, as desired.
\end{IEEEproof}

It is also worth noting that the result of
Lemma~\ref{lemma:ideal:output:pushforward} generalizes to the
nonlinear case under the same arguments, provided that the map
associating the system's initial states to the output samples is
invertible.

\begin{corollary}\longthmtitle{Nonlinear pushforward of output-sample
    reconstructed
    states}\label{corolarry:nonlinear:output:ideal:pushforw}
  Consider the nonlinear system
  \eqref{data:dynamics}-\eqref{data:output} and a sequence of
  trajectories $\xi^i$ as in Assumption~\ref{assumption:sampling},
  Assume that for each $i\in[N^{\flat}:\bar N]$, the map $\mathcal
  H_i:\Phi_{t_i^\ell}(K)\to\Rat{\ell m}$,
  \begin{align*}
    \mathcal H_i(\xi):=(H(t_i^1,\Phi_{t_i^1,t_i^\ell}(\xi)), 
    H(t_i^2,\Phi_{t_i^2,t_i^\ell}(\xi)),\ldots, 
    H(t_i^{\ell},\xi))
  \end{align*}  
  is invertible on its image, where $\Phi_{t_i^\ell}(K)\subset\Rat{d}$
  and $K$ contains the support of the initial conditions'
  distribution. Then, the result of
  Lemma~\ref{lemma:ideal:output:pushforward} remains valid with
  $\bar{\xi}_T^i:=\Phi_{T,t_i^{\ell}}(\mathcal
  H_i^{-1}(\bm{\zeta}^i))$.
\end{corollary}
    
Using Lemma~\ref{lemma:ideal:output:pushforward}, we obtain the
analogue of Proposition~\ref{prop:amb:set:convergence} for
system~\eqref{linear:system} when exact output samples are
assimilated.

\begin{corollary}\longthmtitle{Output-sample based ambiguity radius
    convergence}\label{corollary:output:ambiguity:convergence}
  Under Assumption~\ref{observability:assumption:tv}, further assume
  that the intra-trajectory inter-sampling-time bound $\Delta'$
  satisfies \eqref{Delta:for:sampled:observability} and consider the
  cumulative empirical distribution
  in~\eqref{cumulative:emp:distribution}, where $\bar{\xi}_T^i$ is
  given by \eqref{xiibar:atT:output:dfn} with $N^{\flat}=1$. Then,
  under the hypotheses of Proposition~\ref{prop:amb:set:convergence},
  for any confidence $1-\beta$, the ambiguity radius
  $\eps_N(\beta,\rho_T)$
  satisfies $\lim_{T\to\infty} \eps_N(\beta,\rho_T)=0$.
\end{corollary}

Based on the derived observability results, we also obtain bounds for
the discrepancy between the estimated state from the measurements and
the true state, when the output samples are subject to bounded
observation errors.

\begin{prop}\longthmtitle{State estimation error under bounded
    observation
    errors} \label{prop:state:estimation:output:disturbances}
  Under Assumption~\ref{observability:assumption:tv}, further assume
  that, for each state trajectory $\xi^i$, instead of the exact
  output samples $\bm{\zeta}^i$ in \eqref{xiibar:atT:output:dfn}, we
  measure
  \begin{align}\label{samples:zeta:hat}
    \bm{\widehat{\zeta}}^i=\bm{\zeta}^i+\bm{\delta}^i,
  \end{align}
  for some $\bm{\delta}^i = (\delta_1^i , \dots, \delta_{\ell}^i)$,
  with $\|\delta_l^i\|\le \delta^*$, for $l\in[1:\ell]$.  Also, assume
  that the intra-trajectory inter-sampling-time bound $\Delta'$
  satisfies \eqref{Delta:for:sampled:observability}, or
  \eqref{Delta:for:LTI:sampled:observability} if the system is
  time-invariant.  Then, the estimated state
  $\widehat{\xi}_{t_i^{\ell}}^i:={(\bm W_i\Obs_i)^{\dagger}\bm
    W_i}\widehat{\bm{\zeta}}^i$ satisfies
  \begin{align}\label{xi:difference:bound}
    \|\widehat{\xi}_{t_i^{\ell}}^i-\xi_{t_i^{\ell}}^i\| \le
    \epsilon^*:=
    \sqrt{\frac{\supscr{\tau}{up}}{a\lambda_{\min} 
        \big(\W_{\supscr{\tau}{low}}(t)|_0^{T-\supscr{\tau}{low}}\big)}}\delta^*.
  \end{align}
\end{prop}

The proof of
Proposition~\ref{prop:state:estimation:output:disturbances} is given
in the Appendix. We next provide the analogue to
Theorem~\ref{thm:numerical:ambiguity:radius}, i.e., we determine the
ambiguity radius obtained through the cumulative empirical
distribution $\bar P_{\xi_T}^N$, by pushing forward the estimated
states from the perturbed measurements through the numerical
approximation of the flow. Note that since according to
Assumption~\ref{observability:assumption:tv} the map $t\mapsto A(t)$
is continuous, $\supscr{\Phi}{num}$
satisfies~\eqref{numerical:error:bound} for all $0\le s\le t\le T$,
for some $\fK>0$, and $L:=\max_{t\in[0,T]}\|A(t)\|$.

\begin{thm}\longthmtitle{Ambiguity radius with approximate pushforward
    and observation
    errors}\label{thm:numerical:ambiguity:radius:outputs}
  Under Assumption~\ref{observability:assumption:tv}, let the initial
  condition of~\eqref{linear:system} be supported on the compact
  set~$K$.  Consider a sampling sequence as in
  Assumption~\ref{assumption:sampling}, with inter-trajectory
  sampling-time bound $\Delta>0$, and intra-trajectory
  inter-sampling-time bound $\Delta'>0$. For each effective trajectory
  $\xi^i$, $i\in[N^{\flat}:\bar N]$, we measure the inexact samples
  $\widehat{\bm{\zeta}}^i$ in \eqref{samples:zeta:hat}, and we
  consider the cumulative empirical distribution $\bar P_{\xi_T}^N$
  in~\eqref{cumulative:emp:distribution}, with
  \begin{align*}
    \bar{\xi}_T^i := \supscr{\Phi}{num}_{T,t_i^{\ell}}({(\bm
      W_i\Obs_i)^{\dagger}\bm W_i}\widehat{\bm{\zeta}}^i),
  \end{align*}
  and $\Obs_i$ given by~\eqref{general:observability:matrix:timevar}
  for all $i$. Then, for any confidence $1-\beta$ and $p\ge 1$,
  \eqref{total:ambituity:characterization} holds, with $\psi_N$ as
  in~\eqref{mixed:ambiguity:radius}, $\bar{\eps}_N$ given by
  \begin{align}
    & \bar{\eps}_N(\Delta) := \Big(
    \frac{2^{p-1}}{N}\Big(\frac{(\epsilon^*)^p}{pL\Delta}(e^{pL\Delta
      N}-1)
    \nonumber \\
    & \hspace{8em}+\fK^p\int_1^N(e^{L\Delta
      s}-1)^pds\Big)\Big)^{\frac{1}{p}},
    \label{bar:epsN:output}
  \end{align}
  instead of \eqref{bar:epsN}, and $\epsilon^*$ given by
  \eqref{xi:difference:bound}.
\end{thm}

The proof of this result is given in the Appendix. Under the
hypothesis of Theorem~\ref{thm:numerical:ambiguity:radius:outputs},
the result of Proposition~\ref{prop:effective:horizon} remains valid
with $\bar{\eps}_N$ as given by~\eqref{bar:epsN:output}.
	
\section{Application to UAV Detection}\label{sec:example}

\begin{figure}[tbh]
  \centering
  \subfigure[]{\includegraphics[width=.779\columnwidth]{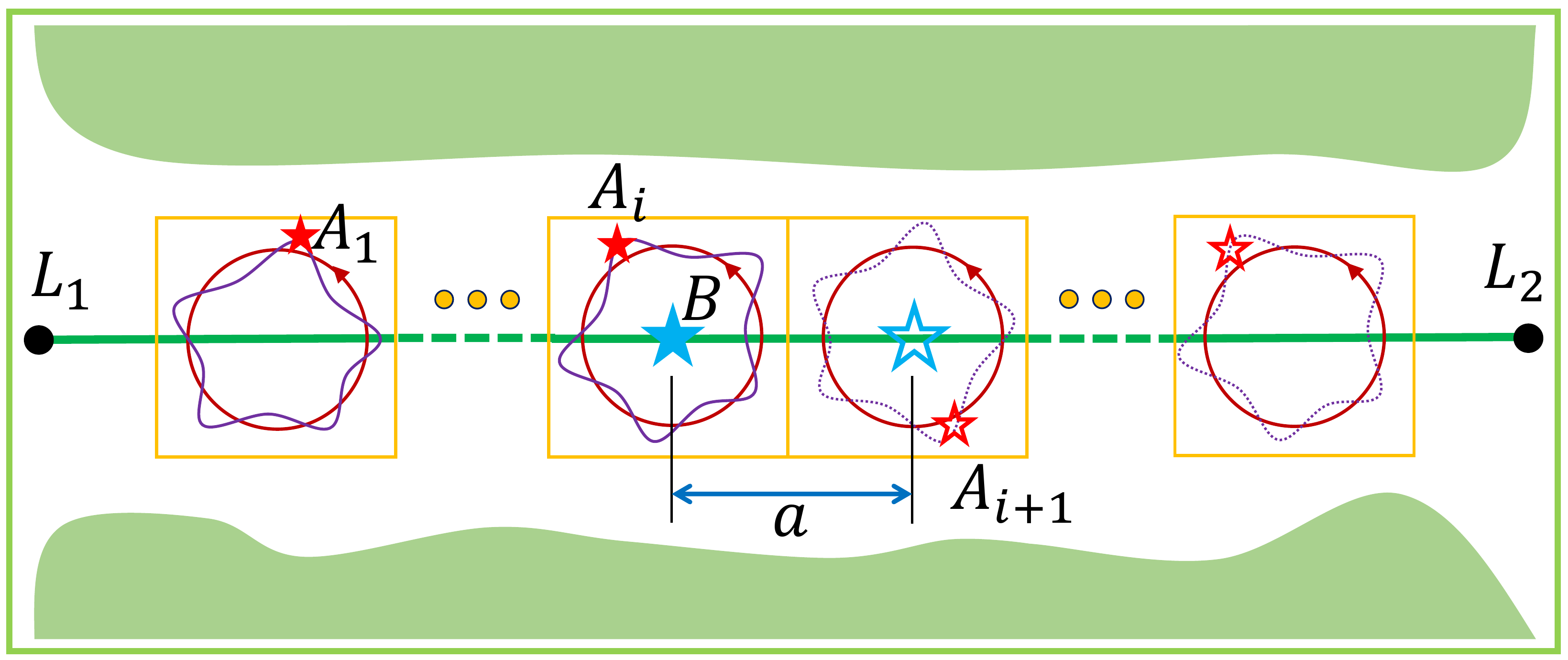}}\hfill
  \subfigure[]{\includegraphics[width=.202\columnwidth]{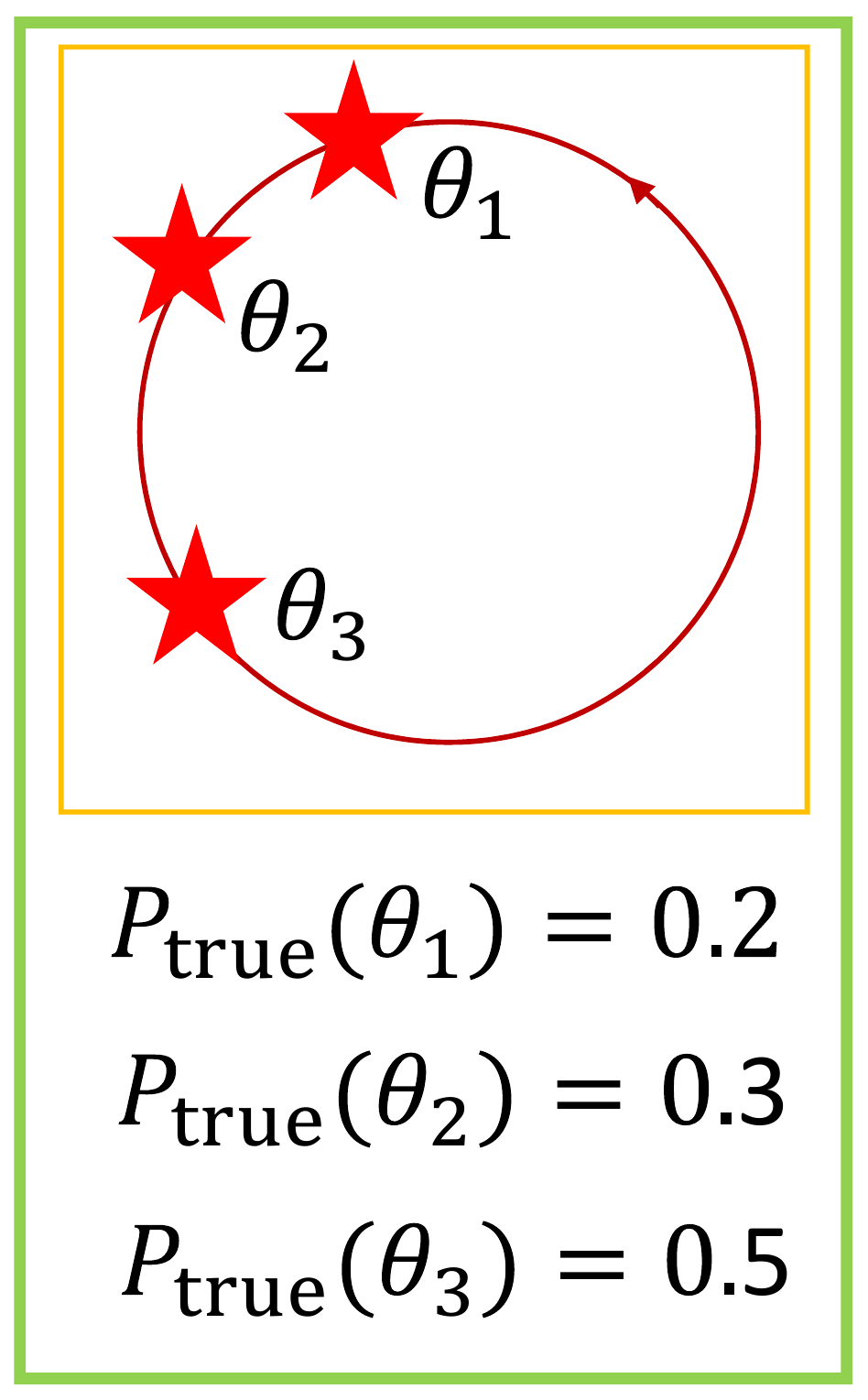}}
  \caption{(a) Red stars depict the location of red UAVs
    $A_1,\ldots,A_i,A_{i+1},\ldots$ at a specific time instant and the
    filled blue star depicts the location of the blue UAV~$B$. Red
    circles represent the trajectories tracked by the red UAVs. The
    non-filled blue star represents the goal position of the blue UAV
    at the end of the optimization horizon, and the non-filled
    stars/dashed trajectories of the red UAVs with index $i+1$ and
    beyond signify that they have not been yet observed by the blue
    UAV. (b) shows the phase angle distribution of the reference
    trajectory tracked by red UAVs, which is unknown to the blue
    UAV.}\label{fig:example}
\end{figure}

Here we illustrate our results in an application scenario involving a
blue UAV that seeks to avoid detection while passing through an area
surveilled by a team of red UAVs, cf. Figure~\ref{fig:example}.
The red team has split the area into squares of identical size and
assigned one UAV per square.  With the vehicle's coordinate system
centered at the middle of the square, each red UAV tracks a circular
orbit $ \chi(t,\theta):=r(\cos(\theta+t),\sin(\theta+t)) $ of radius
$r$ according to the dynamics
\begin{align}\label{pos:vel:dynamics}
  \dot{\xi}_1(t) & =\xi_2(t) \notag
  \\
  \dot{\xi}_2(t) & =\kappa^2({\chi}(t,\theta)-\xi_1(t)).
\end{align}
Here, the phase angle $\theta$ is a parameter not known to
the blue UAV that the red team selects randomly according to some
distribution to make things harder for potential intruders. Note that
this fits into the model~\eqref{data:dynamics} by setting $
\dot{\theta}(t) =0$.  The state of each UAV is
$\xi=(\xi_1,\xi_2,\theta) \in \Rat{5}$.
The purpose of the blue UAV
is to traverse the area along the straight path from location $L_1$ to
location $L_2$ in Figure~\ref{fig:example} while minimizing the
possibility of being detected. The blue UAV is kinematic, actuating
the magnitude of its velocity vector, which points from $L_1$ to $L_2$
and is upper bounded by $v_{\max}$ and lower bounded by $v_{\min}>0$,
since it has not the ability to hover.  The initial conditions of the
red UAVs, i.e., initial position, velocity, and phase angle are
independently sampled from a compactly supported distribution which is
unknown to the blue UAV, except from its support. Thus, vehicle $B$
has no information about the state of each UAV $A_i$ before reaching
region $i$.

\paragraph*{DRO formulation}
While traversing the corresponding square, the blue vehicle collects
(at least three) exact position measurements from $A_i$ before
reaching the center of the square.
At the middle of each square $i$, the blue UAV solves an optimization
problem to tune its velocity up to the center of the next region $i+1$
in order to simultaneously maximize the worst-case distance from UAV
$A_i$, which has already been observed, and the worst-case expected
distance from UAV $A_{i+1}$, based on the collected data of the
preceding UAVs' locations. This maximization is carried out during the
traversal between the centers of squares $i$ and $i+1$, which are at a
distance $a$ apart, and must be performed in $2\pi$ time units (the
surveillance period of the red UAVs). For the purpose of the
optimization, we divide $[0,2\pi]$ into $n$ equal subintervals so
that, at time ${T_i}=i\cdot2\pi$, the blue UAV seeks to determine its
velocity profile $x=(x_1,\ldots,x_n)$ for the next $2\pi$ seconds as
\begin{align*}
  v(T_i+t,x):=x_{\fn}\cdot(1,0), \quad t\in\bigg[(\fn-1) \frac{2\pi}{n},
  \fn \frac{2\pi}{n}\bigg],
\end{align*}
for $\fn=1,\ldots,n$.
%
Here, the velocity profile must belong to
\begin{align*}
  \X:=\{x\in\Rat{n}\,|\, & 0<v_{\min}\le x_{\fn}\le
  v_{\max},\fn=1,\ldots,n, \\
  & \hspace{5.1em} x_1+\cdots+x_n=an/(2\pi)\}.
\end{align*}
Since the blue UAV has no knowledge of the location of UAV $A_{i+1}$,
it uses the previously collected data to solve a DRO formulation that
robustifies its decision against said uncertainty.  Assuming
measurements are exact, it exploits sampled-data observability of the
UAVs' dynamics to recover the states $\xi_{T_i}^1,\ldots,\xi_{T_i}^i$
of the observed red vehicles, and construct an ambiguity set
containing the true distribution of $\xi$ at~$T_i$ with confidence
$1-\beta$. Specifically, using
Corollary~\ref{corolarry:nonlinear:output:ideal:pushforw}, we build an
ambiguity ball $\widehat{\P}_{T_i}^i$ centered at the empirical
distribution $\frac{1}{i}\sum_{k=1}^i\xi_{T_i}^k$ with radius
$\eps_i(\beta)\equiv\eps_i(\beta,\rho_{T_i})$. For each $i=1,2,\ldots$
and associated time $T_i=i\cdot2\pi$, we solve the DRO problem
\begin{align*}
  \sup_{x\in\X}\inf_{P\in \widehat{\P}_{T_i}^i}\bE_P[f_i(x,\xi)],
\end{align*}
where the objective function $f_i$ is given by  
\begin{align*} 
  f_i(x,\xi):=\min_{t\in[T_i,T_i+2\pi]}\bigg\{& \min\Big\{
  \Big\|\xi_{1t}^i-\int_{T_i}^tv(s,x)ds\Big\|^2, \nonumber \\
  &
  \quad\Big\|\supscr{\Phi}{pos+}_t(\xi)-\int_{T_i}^tv(s,x)ds\Big\|^2\Big\}\bigg\},
\end{align*}
with $\xi^i$ the known trajectory of UAV $A_i$,
$\supscr{\Phi}{pos+}_t(\xi):=\subscr{\rm pr}{pos}(\Phi_{\mathit 
t}(\xi))+(a,0)$, and 
$\subscr{\rm
  pr}{pos}(\xi_1,\xi_2,\theta):=\xi_1$.

\paragraph*{Simulation results}
For the simulation results we select the radius $r=1$ for the tracked
trajectories, the square side length $a=2.5$, and set the velocity
bounds for UAV $B$ to $v_{\min}=0.3a/(2\pi)$ and
$v_{\max}=1.5a/(2\pi)$, respectively.  
The tracking (angular) frequency in their
dynamics~\eqref{pos:vel:dynamics} is $\kappa=4$, implying that their
trajectory is periodic with period $2\pi$. Thus, any set of states is
invariant under the flow $\Phi_{T_i}$, for all $T_i=i\cdot2\pi$. The
random values of the phase angle $\theta$ are sampled from the finite
set $\{2.8\pi/4,3.5\pi/4,4.6\pi/4\}$, with the associated
probabilities depicted in Figure~\ref{fig:example}(b), and each UAV is
initiated from the corresponding position
$(\cos(\theta),\sin(\theta))$ with zero velocity, inducing the compact
set of initial states
$K:=\{(\cos(\theta),\sin(\theta),0,0,\theta)\}_{\theta
  \in\{2.8\pi/4,3.5\pi/4,4.6\pi/4\}}\subset\Rat{5}$.
Due to invariance of the flow maps $\Phi_{T_i}$, the corresponding
diameters $\rho_{T_i}$ of the sets $\Phi_{T_i}(K)$ do not change with
$i$. The velocity profile of the blue UAV has $n=4$ subintervals.

\begin{figure}[tbh]
  \centering
  \includegraphics[width=.96\columnwidth]{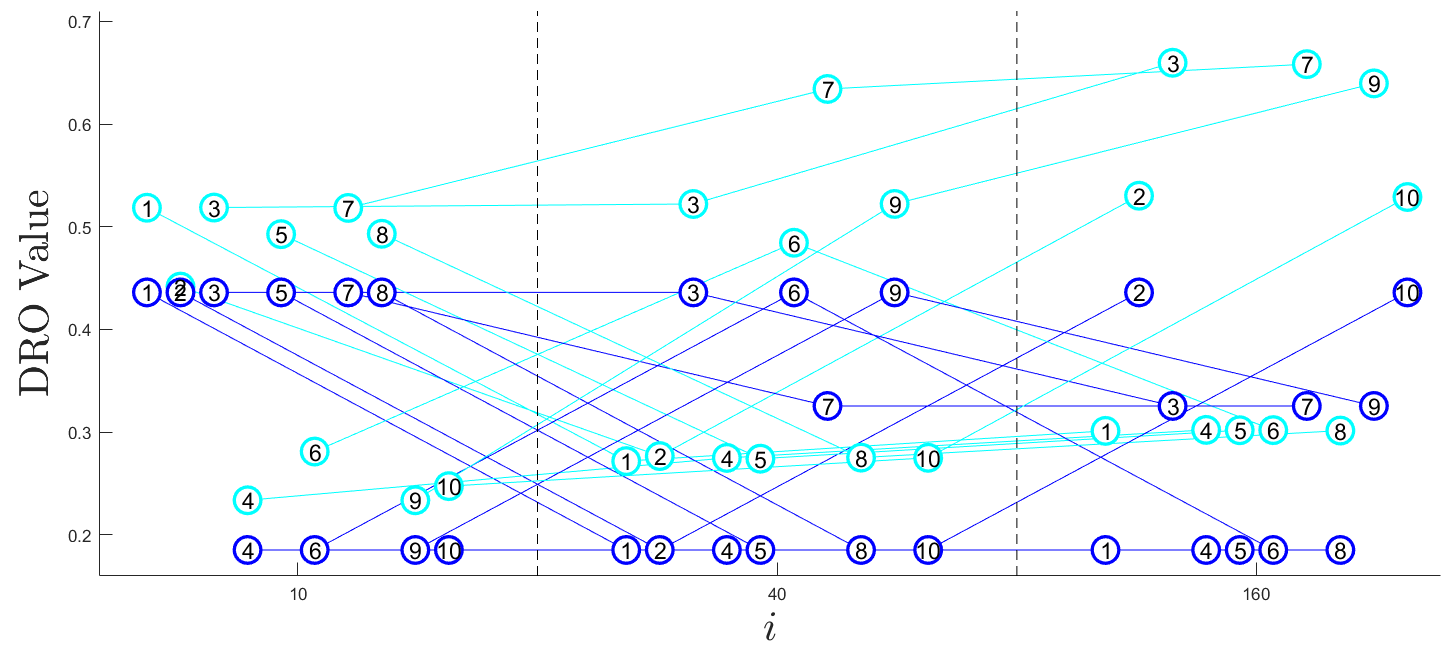}\hfill
  \caption{This plot illustrates the solution of the DRO problem from
    10 independent realizations of the whole detection scenario, when
    the blue UAV crosses the area monitored by the 10th, 40th, and
    160th red UAV, depicted at the left, center, and right region,
    respectively. For each of the 10 realizations the solutions of the
    dynamic and static DRO are illustrated with the cyan and blue
    circles, respectively, and the number of each realization is
    inscribed in the circle. The superiority of the dynamic DRO can be
    seen in the higher optimal values, which statistically increase
    with the number of samples. The fluctuations in the individual
    differences occur due to randomness of the
    realizations.}\label{fig:data:DRO}
\end{figure}

We compare the DRO approach of the paper, where the ambiguity set is
constructed by exploiting all progressively collected samples, with
the static DRO approach, where the ambiguity set is built exclusively
based on the last trajectory's state.  We perform 10 independent
realizations of the detection scenario to illustrate the consistency
of the dynamic DRO benefits throughout each of these independent
experiments.  In each one, we take randomly up to 160 samples from the
probability distribution of $\theta$ to determine the dynamics of the
corresponding red UAVs.  For each realization, we solved the dynamic
DRO using all the samples of the first 10, 40, and 160 red UAVs,
respectively, with the optimal value depicted with cyan in
Figure~\ref{fig:data:DRO}. The static DRO is also solved using only
the single sample of the 10-th, 40-th, and 160-th trajectory (depicted
in blue in Figure~\ref{fig:data:DRO}). It is clear that the
statistical average of the DRO values obtained through the cumulative
empirical distribution outperforms significantly its
single-state-sample counterpart.  Furthermore, the performance of the
suggested DRO scheme improves as the number of samples increases. We
consider a confidence $1-\beta$ such that the ambiguity radius for
$i=10$ trajectory samples is $\eps_i(\beta)=0.17$. From
Corollary~\ref{cor:varying:compact:ambiguity}, the corresponding
radius for $i=40$, $i=160$, and $i=1$ (static DRO), is
$\eps_i(\beta)=0.1201$, $\eps_i(\beta)=0.085$, and
$\eps_i(\beta)=0.3023$, respectively.

\begin{figure}[tbh]
  \centering
  \subfigure[]{\includegraphics[width=.99\columnwidth]{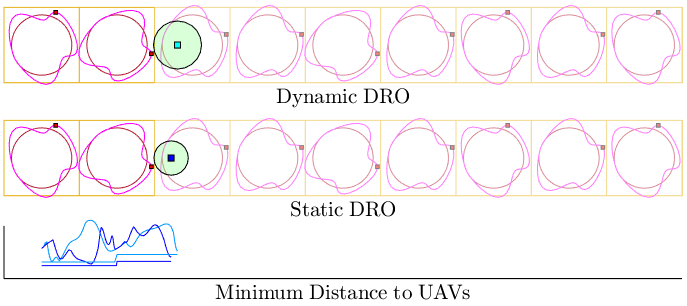}} \\
  \subfigure[]{\includegraphics[width=.49\columnwidth]{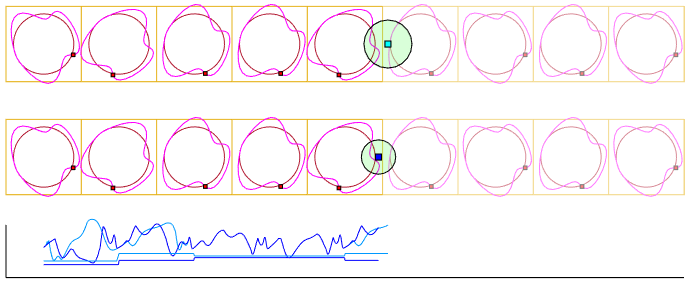}}
  \hfill
  \subfigure[]{\includegraphics[width=.49\columnwidth]{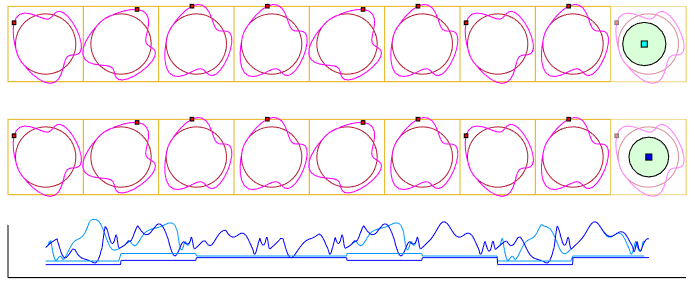}}
  \caption{The figure depicts three
    snapshots of the blue UAV at (a) $T_i=151.6\cdot2\pi$, (b)
    $T_i=154.4\cdot2\pi$, and (c) $T_i=157.6\cdot2\pi$.  The upper and
    lower row of vehicle trajectories in each picture shows the
    outcome of the dynamic and static DRO, respectively. The circles'
    radius is the square root of the DRO value, given also in the
    corresponding piecewise constant plots that are shown below. We
    additionally plotted the evolution of the minimum distance between
    the blue UAV and the red UAVs during the traversal, in cyan and
    blue, for the dynamic and the static DRO, respectively.  The part
    of the snapshots with faded color includes the red UAVs that have
    not been observed by the blue one.} \label{fig:video:instances}
\end{figure}

Figure~\ref{fig:video:instances} shows snapshots from the blue UAV
traversing the area monitored by the 151 to 159th red UAVs in one of
the realizations.  
The radius of the circle moving with vehicle $B$ represents the root of
the minimum squared distance from the red UAVs as obtained from the
solution of the DRO. The dynamic DRO circle is larger than the static
one since a better optimum is obtained in this case due to the
ambiguity set being considerably smaller.

\section{Conclusions}\label{sec:conc}
	
We have developed a framework to compute ambiguity sets of unknown
probability distributions using data collected from dynamically
varying processes that retain the same probabilistic guarantees as
their static counterparts.
Under exact knowledge of the dynamic evolution of the data and
full-state measurements, we have identified conditions on its growth
rate and the sampling rate that ensure the ambiguity sets shrink with
time. In the presence of numerical errors and/or disturbances in the
dynamics, we have also quantified the number of exploitable past
samples necessary to establish the reduction of the ambiguity sets. We
have generalized these results to the case of partial-state
measurements of linear time-varying systems building on their
sample-data observability properties.  Future work will exploit the
construction of time-varying ambiguity sets in receding horizon DRO
problems and explore the extension of the results to consider data
storage limitations, stochastic descriptions of the disturbances,
measurement noise in the observations, and scenarios with incomplete
knowledge of the dynamic evolution of the data.

\appendix 
\section{Appendix}
	
Here we provide proofs of various results of the paper.

\subsection{Proofs from Section~\ref{sec:dyn:ambiguity}}

\begin{IEEEproof}[Proof of Lemma~\ref{lemma:ideal:pushforward}]
  Part (i) follows directly from the fact that
  $\{\xi_0^i\}_{i\in[N^{\flat}:\bar N]}$ are i.i.d. and $\Phi_T$
  is measurable, implying that all $\xi_T^i=\Phi_T\circ\xi_0^i$ are
  also i.i.d. (see \cite[Remark 2.15(iii)]{AK:13}).
  For part (ii), note that
  %
  %
  $\bar P_{\xi_T}^N=\frac{1}{N}\sum_{i=N^{\flat}}^{\bar
    N}\delta_{\bar{\xi}_T^i} =\frac{1}{N}\sum_{i=N^{\flat}}^{\bar
    N}\delta_{\Phi_{T,t_i}(\xi_{t_i}^i)}
  =\frac{1}{N}\sum_{i=N^{\flat}}^{\bar N}\delta_{\xi_T^i}=\widehat
  P_{\xi_T}^N$, as desired, where we have exploited the fact that
  $\Phi_{T,t_i}(\xi_{t_i}^i)=\Phi_{T,t_i} \circ \Phi_{t_i}(\xi_0^i) =
  \Phi_{T}(\xi_0^i)$ in the second to last equality.
  %
\end{IEEEproof}

The proof of Proposition~\ref{prop:dynamics:growth} is based on the
following comparison lemma.
	
\begin{lemma}\longthmtitle{Polynomial growth}\label{lemma:comparison}
  Consider a locally absolutely continuous function
  $\theta:\RgeO\to\RgeO$ which satisfies for almost all $t\ge 0$
  \begin{align}\label{theta:dissipation:assumption}
    \theta(t)>0 \Rightarrow
    \dot{\theta}(t)\le\alpha(t)\theta(t)+M_1\theta(t)^q,
  \end{align}
  for certain $q\in(-\infty,1)$ and $M_1>0$, where the function
  $\alpha$ is locally integrable and
  satisfies~\eqref{alpha:int:hypothesis} for certain $M_2>0$. Then,
  \begin{align}\label{theta:decay:rate}
    \theta(t)\le e^{M_2}(1+M_1(1-q)t)^{\frac{1}{1-q}}(1+\theta_0),
  \end{align}
  where $\theta_0=\theta(0)$.
\end{lemma}
\begin{IEEEproof}
  Let $A(t) =\int_0^t\alpha(s)ds$ and $\bar A(t):=\min_{0\le s\le
    t}A(s)$.  Since $A(0)=0$, also $\bar A(0)=0$. Furthermore, $\bar
  A$ is locally absolutely continuous due to the fact that the same
  property holds for $A$. From the latter, combined with the fact that
  $\bar A$ is nonincreasing and satisfies $\bar A(0)=0$, it follows
  that there exists a nonnegative locally integrable function
  $\bar{\alpha}$ with $\bar A(t)=-\int_0^t\bar{\alpha}(s)ds$, for all
  $t\ge 0$.  Then, if we define
  $\widehat{\alpha}(t):=\alpha(t)+\bar{\alpha}(t)$, $t\ge 0$, we
  obtain for all $t\ge 0$,
  \begin{subequations}
    \begin{align}
      \widehat{\alpha}(t) & \ge \alpha(t),
      \label{alpha:hat:vs:alpa}	
      \\
      \int_0^t\widehat{\alpha}(s) ds & \ge 0,
      \label{alpha:hat:int:lowbound}
      \\
      \int_0^t\widehat{\alpha}(s) ds & \le M_2,
      \label{alpha:hat:int:upbound}
    \end{align}
  \end{subequations}
  Indeed,~\eqref{alpha:hat:vs:alpa} follows directly from the fact
  that $\bar{\alpha}$ is nonnegative. In addition, we have that
  $\int_0^t\widehat{\alpha}(s)ds=\int_0^t\alpha(s)ds+\int_0^t\bar{\alpha}(s)ds=A(t)-\bar
  A(t)=A(t)-\min_{0\le s\le t}A(s)\ge 0$, which establishes
  \eqref{alpha:hat:int:lowbound}. Finally, given that $\bar
  A(t)=\min_{0\le s\le t}A(s)=\int_0^{\tau} \alpha(s)ds$ for some
  $\tau\in[0,t]$, we get that
  $\int_0^t\widehat{\alpha}(s)ds=\int_0^t\alpha(s)ds-\bar
  A(t)=\int_0^t\alpha(s)ds-\int_0^{\tau} \alpha(s)ds=\int_{\tau}^t
  \alpha(s)ds\le M_2$, because of \eqref{alpha:int:hypothesis}.
  From~\eqref{theta:dissipation:assumption} and
  \eqref{alpha:hat:vs:alpa}, we obtain for almost all $t\ge 0$,
  $\theta(t)>0 \Rightarrow \dot{\theta}(t)
  \le\widehat{\alpha}(t)\theta(t)+M_1\theta(t)^q$.  Hence, by defining
  \begin{align}\label{sigma:dfn}
    \sigma(t):=\theta(t/M_1), t\ge 0,
  \end{align}
  it follows that for almost all $t\ge 0$,
  \begin{align}\label{sigma:diff:ineq}
    \sigma(t)>0 \Rightarrow \dot{\sigma}(t) =
    \frac{1}{M_1}\dot{\theta}(t/M_1) \le
    \gamma(t)\sigma(t)+\sigma(t)^q,
  \end{align}
  with $\gamma(t):=1/M_1\widehat{\alpha}(t/M_1)$. Then, we get from
  \eqref{alpha:hat:int:lowbound} that
  \begin{align}\label{gamma:low:bound}
    \int_0^t\gamma(s)ds\ge 0,\quad\forall t\ge 0,
  \end{align}
  and from \eqref{alpha:hat:int:upbound} that
  \begin{align}\label{gamma:up:bound}
    \int_0^t\gamma(s)ds = \int_0^t\frac{1}{M_1} \hspace*{-.5ex}
    \widehat{\alpha}\Big(\frac{s}{M_1}\Big)ds = \int_0^{t/M_1}
    \hspace*{-2.5ex} \widehat{\alpha}(\tau)d\tau\le M_2,
  \end{align}
  for all $t\ge 0$. Now, let $\lambda$ be the solution of 
  \begin{align}\label{lambda:diff:ineq}
    \dot{\lambda}(t)=\gamma(t)\lambda(t)+\lambda(t)^q,\quad 
    \lambda(0)=\sigma(0)+1,
  \end{align}
  which is defined for all $t\ge 0$ and is nondecreasing. We claim
  \begin{align}\label{lambda:vs:sigma}
    \lambda(t)\ge\sigma(t),\quad\textup{for all}\; t\ge 0.
  \end{align}
  Indeed, suppose on the contrary that $\sigma(T)>\lambda(T)$ for certain 
  $T>0$. Define $\tau:=\inf\{\bar t\ge 0\,|\,\sigma(t)>\lambda(t),\forall 
  t\in(\bar t,T]\}$. Then,
  \begin{subequations}
    \begin{align}
      \sigma(t) & >\lambda(t),\;\forall t\in(\tau,T],\label{theta:vs:lambda:after:tau} 
      \\
      \sigma(\tau) & =\lambda(\tau).\label{theta:vs:lambda:at:tau}
    \end{align}
  \end{subequations}
  By \eqref{sigma:diff:ineq}, \eqref{lambda:diff:ineq},
  \eqref{theta:vs:lambda:at:tau}, and the comparison lemma~\cite[Lemma
  3.4]{HKK:02}, we get that $\sigma(t)\le\lambda(t)$ for all
  $t\in(\tau,T]$, contradicting \eqref{theta:vs:lambda:after:tau}. We
  next show that the absolutely continuous function
  \begin{align}\label{eta:dfn}
  \eta(t):=e^{\int_0^t\gamma(s)ds}(1+t(1-q))^{\frac{1}{1-q}}(1+\theta_0),t\ge 0,
  \end{align}
  satisfies
  \begin{align}\label{eta:diff:ineq}
    \dot{\eta}(t)\ge\gamma(t)\eta(t)+\eta(t)^q
  \end{align} 
  for almost all $t\ge 0$. Indeed, note that $\dot{\eta}(t) =
  \gamma(t)\eta(t)+e^{\int_0^t\gamma(s)ds}(1+\theta_0)(1+t(1-q))^{\frac{q}{1-q}}$
   for almost all $t\ge 0$. Solving with respect to $1+t(1-q)$ in 
  \eqref{eta:dfn}, we deduce
  \begin{align*}
    \dot{\eta}(t) & =\gamma(t)\eta(t)+
    e^{\int_0^t\gamma(s)ds}(1+\theta_0)\bigg(\frac{\eta(t)}{e^{\int_0^t\gamma(s)ds}(1
      + \theta_0)}\bigg)^q \\ 
    & 
    =\gamma(t)\eta(t)+\big(e^{\int_0^t\gamma(s)ds}(1+\theta_0)\big)^{1-q}\eta(t)^q,
  \end{align*}	
  which implies \eqref{eta:diff:ineq}, due to \eqref{gamma:low:bound}
  and the fact that $1-q>0$.  Thus, since by \eqref{sigma:dfn}
  and~\eqref{lambda:diff:ineq} it holds that
  $\eta(0)=1+\theta(0)=1+\sigma(0)=\lambda(0)$, we get from
  \eqref{lambda:diff:ineq}, \eqref{eta:diff:ineq}, and the comparison
  lemma~\cite[Lemma 3.4]{HKK:02}, that $\eta(t)\ge\lambda(t)$ for all
  $t\ge 0$. Hence, from \eqref{lambda:vs:sigma} and the definition of
  $\sigma$, we have that $\eta(t)\ge \theta(\frac{t}{M_1})$ for all
  $t\ge 0$. From the latter, \eqref{gamma:up:bound}, and
  \eqref{eta:dfn}, we deduce \eqref{theta:decay:rate}, which completes
  the proof.
\end{IEEEproof}


We are ready to prove Proposition~\ref{prop:dynamics:growth}.
		
\begin{IEEEproof}[Proof of Proposition~\ref{prop:dynamics:growth}]
  Consider an initial condition $\xi_0\in\Rat{d}$ and let $\xi$ denote
  the solution of \eqref{data:dynamics} defined for all $t\ge
  0$.
  Then, $t\mapsto V(\xi(t))$ is differentiable, and from
  \eqref{Lyapunov:inequality},
  \begin{align*}
    V(\xi(t))>0 \Rightarrow \frac{d}{dt}V(\xi(t)) &
    =DV(\xi(t))F(t,\xi(t))
    \\
    & \le \alpha(t)V(\xi(t))+M_1V(\xi(t))^q.
  \end{align*}
  Then, if $M_1=0$, we obtain by the comparison lemma \cite[Lemma
  3.4]{HKK:02}
  that $V(\xi(t))\le V(\xi_0)e^{\int_0^t\alpha(s)ds}$ for all
  {$t\ge 0$}, and thus, by \eqref{sandwitch:bounds}, that
  \eqref{decay:M1:eq:zero} holds. If $M_1>0$, then, by using
  Lemma~\ref{lemma:comparison} with $\theta(t):=V(\xi(t))$, we have
  from \eqref{theta:decay:rate} that
  $V(\xi(t))\le e^{M_2}(1+M_1(1-q)t)^{\frac{1}{1-q}}(1+V(\xi_0))$ 
  for all {$t\ge 0$}.
  %
  Consequently,
  {it follows} from \eqref{sandwitch:bounds} that
  \begin{align*}
    \|\xi(t)\|\le(e^{M_2}(1+a_2\|\xi_0\|^r)/a_1)^{\frac{1}{r}}(1+M_1(1-q)t)^{\frac{1}{r(1-q)}}
  \end{align*}
  for all $t\ge 0$, which establishes \eqref{xi:growth:rate} with
  $\bar M$ and $\bar c$ as given by \eqref{constant:K}. The proof is
  complete.
\end{IEEEproof}

\begin{IEEEproof}[Proof of Fact II in
  Proposition~\ref{prop:amb:set:convergence}]
  It suffices to show that $\lim_{x\to 0}h^{-1}(ax)/x^{\frac{1}{\bar
      q}}=0$, or equivalently, that
  \begin{align*}
    h^{-1}(ax)^\frac{\bar q+2}{2}/(x^{\frac{1}{\bar q}})^\frac{\bar
      q+2}{2}\to 0.
  \end{align*}
  Using L'H\^opital's rule and $\bar q>2$,
  it follows that there exists $\bar y>0$ with $h(y)>y^\frac{\bar q+2}{2}$ for 
  all $0<y\le\bar y$. Thus,
  \begin{align*}
    h^{-1}(ax)^\frac{\bar q+2}{2}/(x^{\frac{1}{\bar q}})^\frac{\bar
      q+2}{2} & <h(h^{-1}(ax))/(x^{\frac{1}{\bar q}})^\frac{\bar
      q+2}{2}
    \\
    & =ax/x^{\frac{\bar q+2}{2\bar q}}=ax^{\frac{\bar q-2}{2\bar
        q}}\to 0,
  \end{align*}
  since $\bar q>2$, and we get the result.
\end{IEEEproof}

\subsection{Proofs from
  Section~\ref{sec:errors:disturbances}}

The proof of Theorem~\ref{thm:numerical:ambiguity:radius} relies on
the following elementary result on the Wasserstein distance between
two discrete distributions with the same number of elements.

\begin{lemma}\longthmtitle{Wasserstein distance of discrete
    distributions}\label{lemma:two:empirical:Wdist} 
  Consider the finite sequences $(X_i)_{i=1}^N$, $(Y_i)_{i=1}^N$ in
  $\Rat{d}$ and the corresponding discrete distributions
  $\widehat{\mu}_X^N=\frac{1}{N}\sum_{i=1}^N\delta_{X_i}$,
  $\widehat{\mu}_Y^N=\frac{1}{N}\sum_{i=1}^N\delta_{Y_i}$. Then, it
  holds that $W_p(\widehat{\mu}_X^N,\widehat{\mu}_Y^N)\le
  (\frac{1}{N}\sum_{i=1}^N\|X_i-Y_i\|^p)^\frac{1}{p}$.
\end{lemma}
%
%
\begin{IEEEproof}
  Consider the probability measure
  $\pi=\frac{1}{N}\sum_{i=1}^N\delta_{(X_i,Y_i)}$ on
  $\Rat{d}\times\Rat{d}$.  Then, it follows that $\widehat{\mu}_X^N$,
  $\widehat{\mu}_Y^N$ are marginals of $\pi$ and that
  $W_p(\widehat{\mu}_X^N,\widehat{\mu}_Y^N)
  \le(\int_{\Rat{d}\times\Rat{d}}\|x-y\|^p\pi(dx,dy))^\frac{1}{p}
  =(\frac{1}{N}\sum_{i=1}^N\|X_i-Y_i\|^p)^\frac{1}{p}$,
  as claimed.
\end{IEEEproof}

\begin{IEEEproof}[Proof of
  Theorem~\ref{thm:numerical:ambiguity:radius}]
  Notice that the true distribution $P_{\xi_T}$ of the system state at
  $T$ will be supported in the compact set $\Phi_T(K)$. Thus, we get
  from \eqref{rho:T} and Corollary~\ref{cor:varying:compact:ambiguity}
  that for the given confidence $1-\beta$, the Wasserstein distance
  between the ideal empirical distribution $\widehat P_{\xi_T}^N$ in
  \eqref{emp:distribution:atT} and $P_{\xi_T}$ will satisfy
  $\bP(W_p(\widehat P_{\xi_T}^N,P_{\xi_T})\le\eps_N(\beta,\rho_T))\ge
  1-\beta$. Based on the latter, to establish
  \eqref{total:ambituity:characterization}, it suffices to show that
  \begin{align} 
    W_p(\bar P_{\xi_T}^N,\widehat P_{\xi_T}^N) & \le
    \bar{\eps}_N(\Delta),
    \label{distace:between:hatP:barP}
  \end{align}
  and take into account \eqref{mixed:ambiguity:radius} and the
  triangle inequality for $W_p$.
  
  To show \eqref{distace:between:hatP:barP}, recall that $\widehat
  P_{\xi_T}^N=\frac{1}{N}\sum_{i=N^{\flat}}^{\bar N}\delta_{\xi_T^i}$
  and that $\xi_T^i=\Phi_{T,t_i}(\xi_{t_i}^i)$ for each $i$. Then,
  given that $t_{i+1}-t_i\le\Delta$ for each $i$ and taking into
  account \eqref{numerical:error:bound}, and that $t_{\bar N}=T$, we
  get from Lemma~\ref{lemma:two:empirical:Wdist} that
  \begin{align*}
    & W_p(\bar P_{\xi_T}^N,\widehat P_{\xi_T}^N)
    \le\Big(\frac{1}{N}\sum_{i=N^{\flat}}^{\bar N}\|\xi_T^i-\bar
    \xi_T^i\|^p\Big)^\frac{1}{p}
    \\
    & =\Big(\frac{1}{N}\sum_{i=N^{\flat}}^{\bar N}
    \|\Phi_{T,t_i}(\xi_{t_i}^i)-\supscr{\Phi}{num}_{T,t_i}(\xi_{t_i}^i)\|^p\Big)^\frac{1}{p}
    \\
    & \le \Big(\frac{1}{N}\sum_{i=N^{\flat}}^{\bar N}
    \fK^p(e^{L(T-t_i)}-1)^p\Big)^\frac{1}{p}
    \\
    & \le \fK
    \Big(\frac{1}{N}\sum_{i=1}^N(e^{L\Delta(N-i)}-1)^p\Big)^\frac{1}{p}
    =\fK\Big(\frac{1}{N}\sum_{i=1}^{N-1}(e^{L\Delta
      i}-1)^p\Big)^\frac{1}{p}.
  \end{align*}
  Then, the result is a consequence of the derived bound and the fact that for 
  any $a>0$, $p\ge 1$, and $N\in\bN$, it holds that $\sum_{i=1}^{N-1}(e^{ai}-1)^p 
  \le\int_1^N(e^{as}-1)^p ds$.
\end{IEEEproof}

\subsection{Proofs from
  Section~\ref{sec:partial:observations}}

\begin{IEEEproof}[Proof of
  Lemma~\ref{lemma:sampled:observability:LTI}]
  Under \textbf{H1} in
  Assumption~\ref{sampled:observability:assumptions},
  \begin{align}\label{Osample:eq:Operiodic}
    \Obs_ie^{A(t_i^{\ell}-t_i^1)}
    =\Obs_{\ell}(e^{A\Delta'},C):=\left(\begin{matrix}
        C \\
        Ce^{A\Delta'} \\
        \vdots \\
        C\big(e^{A\Delta'}\big)^{\ell-1}
      \end{matrix}\right).
  \end{align}
  In addition, due to \textbf{H1}, it follows that system
  \eqref{linear:system} is $\Delta'$-sampled observable
  \cite[Proposition~6.2.11]{EDS:98}, namely, the pair
  $(e^{A\Delta'},C)$ is observable.  Thus, the observability matrix
  $\Obs_d(e^{A\Delta'},C)$ corresponding to this pair has full rank
  $d$. Since $\ell\ge d$, $\Obs_{\ell}(e^{A\Delta'},C)$ is also of
  full rank $d$, which by~\eqref{Osample:eq:Operiodic}, establishes
  left invertibility of $\Obs_i$.  When \textbf{H2} holds, the result
  follows analogously from
  the Corollary in Section IV of \cite{GK:99}, which
  implies that
  \begin{align*}
    {\rm rank}(\Obs_i)\equiv{\rm rank}\Big(\Obs_{t_i^1\cdots
      t_i^{(\bar d+1) d}}^{\exp-}\Big)=d.
  \end{align*} 
  Thus, since $\ell\ge (\bar d+1)d$, we also obtain that $\Obs_i$ is
  of full rank, and hence, left invertible. Finally, the same approach
  can be followed under \textbf{H3}, with left invertibility of
  $\Obs_i$ guaranteed by~\cite[Theorem~2]{SZ-FA:16}.
\end{IEEEproof}

Before proceeding with the proof of
Proposition~\ref{prop:robust:observability}, we state some
intermediate results to ensure robust invertibility of the
sample-observability matrix $\Obs_i$. Given $i\in[1:\bar N]$, let
\begin{align}\label{integrals:Wtaul}
  {\W_{\tau(i)}(t_i^1)\big|_{[t_i^l,t_i^{l+1}]}} & :=
  \int_{t_i^l}^{t_i^l+\tau_l(i)}K(s,t_i^1+\tau(i))ds,
\end{align}
for each $l\in[1:\ell-1]$, with the intra-trajectory
inter-sampling-time lengths $\tau_l(i)=t_i^{l+1}-t_i^l$ as
in~\eqref{times:hatt:tau}, and $\tau(i)=t_i^\ell-t_i^1$. Also, recall
that $K_s(s,t)=\frac{\partial}{\partial s}K(s,t)$, which is continuous
because of the regularity hypotheses of
Assumption~\ref{observability:assumption:tv}. 

\begin{lemma}\longthmtitle{Integral-limits
    inequalities}\label{lemma:integral:limits:ineq}
  For each $l\in[1:\ell-1]$, the integral
  ${\W_{\tau(i)}(t_i^1)\big|_{[t_i^l,t_i^{l+1}]}}$ satisfies
  \begin{align*}
    \Big\|\W_{\tau(i)}(t_i^1)\big|_{[t_i^l,t_i^{l+1}]}
    -\frac{\tau_l(i)}{2}( & K(t_i^l,t_i^1+\tau(i)) \\
    & +K(t_i^{l+1},t_i^1+\tau(i)))\Big\| \\
    \le &
    \frac{\tau_l^2(i)}{4}\max_{s\in[t_i^l,t_i^{l+1}]}\|K_s(s,t_i^1+\tau(i))\|.
  \end{align*} 
\end{lemma}
\begin{IEEEproof}
  From the mean value inequality,
  \begin{align*}
    & \Big\|\W_{\tau(i)}(t_i^1)\big|_{[t_i^l,t_i^{l+1}]}
    -\frac{\tau_l(i)}{2}(K(t_i^l,t_i^1+\tau(i)) \\
    & \hspace{11.2em}+K(t_i^{l+1},t_i^1+\tau(i)))\Big\|  \\
    &
    =\bigg\|\int_{t_i^l}^{t_i^l+\frac{\tau_l(i)}{2}}(K(s,t_i^1+\tau(i))
    -K(t_i^l,t_i^1+\tau(i)))ds \\
    &
    \quad+\int_{t_i^l+\frac{\tau_l(i)}{2}}^{t_i^l+\tau_l(i)}(K(s,t_i^1+\tau(i))-
    K(t_i^{l+1},t_i^1+\tau(i)))ds\bigg\| \\
    & \le \int_{t_i^l}^{t_i^l+\frac{\tau_l(i)}{2}}(s-t_i^l)ds
    \max_{s\in[t_i^l,t_i^l+\frac{\tau_l(i)}{2}]}\|K_s(s,t_i^1+\tau(i))\| \\
    & \quad
    +\int_{t_i^l+\frac{\tau_l(i)}{2}}^{t_i^l+\tau_l(i)}\Big(s-t_i^l
    -\frac{\tau_l(i)}{2}\Big)ds \\
    &
    \hspace{8.2em}\times\max_{s\in[t_i^l+\frac{\tau_l(i)}{2},t_i^{l+1}]}
    \|K_s(s,t_i^1+\tau(i))\| \\
    & \le
    \frac{\tau_l^2(i)}{4}\max_{s\in[t_i^l,t_i^{l+1}]}\|K_s(s,t_i^1+\tau(i))\|,
  \end{align*} 
  %
  establishing the result.
\end{IEEEproof}

The following result provides an upper bound for the distance between
the matrices $\Obs_i^{\top}\bm W_i\bm W_i \Obs_i$ and
$\W_{\tau(i)}(t_i^1)$ in the induced Euclidean norm.

\begin{lemma}\longthmtitle{Distance between $\Obs_i^{\top}\bm W_i\bm
    W_i \Obs_i$ and
    $\W_{\tau(i)}(t_i^1)$}\label{lemma:distance:to:Gramian}
  Let $\W_{\tau(i)} (t_i^1)$ and $\bm W_i$ as given by \eqref{Wtau} and
  \eqref{weight:matrix:W}, and consider the sample-observability
  matrix $\Obs_i$ in \eqref{general:observability:matrix:timevar} and
  an intra-trajectory inter-sampling-time bound $\Delta'>0$. Then,
  \begin{align}\label{Garmian:closeto:Obs}
    \|\Obs_i^{\top}\bm W_i\bm W_i\Obs_i & -\W_{\tau(i)}(t_i^1)\|
    \notag
    \\
    \le & \frac{\tau(i)\Delta'}{4}
    \max_{s\in[t_i^1,t_i^1+\tau(i)]}\|K_s(s,t_i^1+\tau(i))\|.
  \end{align}
\end{lemma}
\begin{IEEEproof}
  From the {definitions of $\Obs_i$ in
    \eqref{general:observability:matrix:timevar}, $K(s,t)$ in
    \eqref{map:K}, $\W_{\tau(i)}(t_i^1)$ in \eqref{Wtau}, the
    $\tau_l(i)$'s in \eqref{times:hatt:tau}, the matrix $\bm W_i$ in
    \eqref{weight:matrix:W}, and the matrices
    $\W_{\tau(i)}(t_i^1)\big|_{[t_i^l,t_i^{l+1}]}$ in
    \eqref{integrals:Wtaul},}
  we have 
  \begin{align}
    & \|\Obs_i^{\top}\bm W_i\bm W_i\Obs_i-\W_{\tau(i)}(t_i^1)\|
    \nonumber
    \\
    & = \|\left(\begin{matrix}
        w_1(i)\Phi(t_i^1,t_i^{\ell})^{\top}C(t_i^1)^{\top} & \cdots &
        w_{\ell}(i)\Phi(t_i^{\ell},t_i^{\ell})^{\top}C(t_i^{\ell})^{\top}
      \end{matrix}\right) \nonumber 
    \\
    & \hspace{8.5em}\left.\times\left(\begin{matrix}
          w_1(i)C(t_i^1)\Phi(t_i^1,t_i^{\ell})\\
          \vdots \\
          w_{\ell}(i)C(t_i^{\ell})\Phi(t_i^{\ell},t_i^{\ell})
        \end{matrix}\right)-\W_{\tau(i)}(t_i^1)\right\| \nonumber 
    \\
    &
    =\bigg\|\sum_{l=1}^{\ell}w_l^2(i)K(t_i^l,t_i^1+\tau(i))-\W_{\tau(i)}(t_i^1)\bigg\|
    \nonumber
    \\
    &
    =\bigg\|\sum_{l=1}^{\ell-1}\frac{\tau_l(i)}{2}(K(t_i^l,t_i^1+\tau(i))
    +K(t_i^{l+1},t_i^1+\tau(i))) \nonumber
    \\
    & \hspace{18em} -\W_{\tau(i)}(t_i^1)\bigg\|\nonumber
    \\
    & \le
    \sum_{l=1}^{\ell-1}\Big\|\frac{\tau_l(i)}{2}(K(t_i^l,t_i^1+\tau(i))+K(t_i^{l+1},t_i^1
    +\tau(i))) \nonumber
    \\
    & \hspace{16em}-\W_{\tau(i)}(t_i^1)\big|_{[t_i^l,t_i^{l+1}]}\Big\|
    \nonumber
    \\
    &\le\sum_{l=1}^{\ell-1}\frac{\tau_l^2(i)}{4}
    \max_{s\in[t_i^l,t_i^{l+1}]}\|K_s(s,t_i^1+\tau(i))\|,\label{Gramian:vs:Obs}
  \end{align}
  where we have used Lemma~\ref{lemma:integral:limits:ineq} in the
  last inequality.  Thus, from the bound on the maximum inter-sampling
  time, we get \eqref{Garmian:closeto:Obs}, which establishes the
  result.
\end{IEEEproof}

When the system is time invariant, i.e., $A(t)\equiv A$ and
$C(t)\equiv C$, we obtain the following corollary to
Lemma~\ref{lemma:distance:to:Gramian}, with a more explicit bound for
$\|\Obs_i^{\top}\bm W_i\bm W_i \Obs_i-\W_{\tau(i)}(t_i^1)\|$.

\begin{corollary}\longthmtitle{Distance between $\Obs_i^{\top}\bm
    W_i\bm W_i\Obs_i$ and $\W_{\tau(i)}(t_i^1)$ for LTI
    systems}\label{corollary:distance:to:Gramian:LTI}
  Under the assumptions of Lemma~\ref{lemma:distance:to:Gramian},
  when~\eqref{linear:system} is time invariant,
  \begin{align*} 
    \|\Obs_i^{\top}\bm W_i\bm W_i\Obs_i- & \W_{\tau(i)}(t_i^1)\|
    \\
                                         & \le \frac{\tau(i)\Delta'}{2}\max_{s\in[0,\tau(i)]}\|\widehat
                                           K(s-\tau(i))A\|.
  \end{align*}
\end{corollary}
\begin{IEEEproof}
  From~\eqref{function:hatK}, note that $\dot{\widehat K}(s)=
  \frac{d}{ds}(e^{A^{\top}s}C^{\top}Ce^{As})
  =A^{\top}e^{A^{\top}s}C^{\top}Ce^{As}+e^{A^{\top}s}C^{\top}Ce^{As}A$.
  Since the spectral norm of a matrix equals that of its transpose, it
  follows that $\|\dot{\widehat K}(s)\|\le 2\|\widehat K(s)A\|$.
  %
  %
  %
  %
  Combining this with the bound for $\|\Obs_i^{\top}\bm W_i\bm W_i
  \Obs_i-\W_{\tau(i)}(t_i^1)\|$ in
  Lemma~\ref{lemma:distance:to:Gramian}, and the fact that
  $K(s,t)=\widehat K(s-t)$,
  \begin{align*}
    &
    \frac{\tau(i)\Delta'}{4}\max_{s\in[t_i^1,t_i^1+\tau(i)]}\|K_s(s,t_i^1+\tau(i))\|
    \\
    & =
    \frac{\tau(i)\Delta'}{4}\max_{s\in[t_i^1,t_i^1+\tau(i)]}\|\dot{\widehat
      K}(s-(t_i^1+\tau(i)))\|
    \\
    & = \frac{\tau(i)\Delta'}{4}\max_{s\in[0,\tau(i)]}\|\dot{\widehat
      K}(s-\tau(i))\|
    \\
    & \le \frac{\tau(i)\Delta'}{2}\max_{s\in[0,\tau(i)]}\|\widehat
    K(s-\tau(i))A\|,
  \end{align*}
  which completes the proof.  
\end{IEEEproof}


Based on the obtained results, we prove 
Proposition~\ref{prop:robust:observability}.  

\begin{IEEEproof}[Proof of
  Proposition~\ref{prop:robust:observability}]
  Let $i\in[1:\bar N]$ and note that for any pair of symmetric
  matrices $P$ and $Q$,
  \begin{align*}
    \lambda_{\min}(P) & =\min_{\|x\|=1}x^{\top}Px=\min_{\|x\|=1}x^{\top}(Q+P-Q)x \\
    & \le \min_{\|x\|=1}x^{\top}Qx+\max_{\|x\|=1}x^{\top}(P-Q)x \\
    & =\lambda_{\min}(Q) +\|P-Q\|.
  \end{align*}
  Combining this observation with the result of
  Lemma~\ref{lemma:distance:to:Gramian},
  \begin{align*}
    \lambda_{\min}(\W_{\tau(i)}(t_i^1))\le & \;
    \lambda_{\min}(\Obs_i^{\top}\bm
    W_i\bm W_i\Obs_i) \\
    & +\|\Obs_i^{\top}\bm W_i\bm W_i\Obs_i-\W_{\tau(i)}(t_i^1)\|\\
    \le & \; \lambda_{\min}(\Obs_i^{\top}\bm W_i\bm W_i\Obs_i) \\
    &+\frac{\tau(i)\Delta'}{4}\max_{s\in[t_i^1,t_i^1+\tau(i)]}\|K_s(s,t_i^1+\tau(i))\|.
  \end{align*}
  Thus, we get from \eqref{minimum:eigenvalue:over:Gramians} and 
  by plugging in \eqref{Delta:for:sampled:observability}  that 
  \begin{align*}
    & \lambda_{\min}(\Obs_i^{\top}\bm W_i\bm W_i\Obs_i) \\
    & \ge 
    \lambda_{\min}(\W_{\tau(i)}(t_i^1))-\frac{\tau(i)\Delta'}{4}\max_{s\in[t_i^1,t_i^1+\tau(i)]}
    \|K_s(s,t_i^1+\tau(i))\| \\
    & \ge \lambda_{\min}(\W_{\tau(i)}(t_i^1))- 
    \frac{\tau(i)(1-a)\lambda_{\min}\big(\W_{\supscr{\tau}{low}} 
      (t)|_0^{T-\supscr{\tau}{low}}\big)}{\supscr{\tau}{up}} \\
    & \hspace{7.5em}\times\frac{\max_{s\in[t_i^1,t_i^1+\tau(i)]} 
      \|K_s(s,t_i^1+\tau(i))\|}{\max_{\supscr{\tau}{low}\le t\le 
        T,\max\{0,t-\supscr{\tau}{up}\}\le s\le t}\|K_s(s,t)\|} \\
    & \ge \lambda_{\min}(\W_{\tau(i)}(t_i^1))- 
    \frac{\tau(i)}{\supscr{\tau}{up}}(1-a) 
    \lambda_{\min}\big(\W_{\supscr{\tau}{low}}(t)|_0^{T-\supscr{\tau}{low}}\big) \\
    & \ge\lambda_{\min}(\W_{\tau(i)}(t_i^1))-(1-a)\lambda_{\min} 
    \big(\W_{\supscr{\tau}{low}}(t)|_0^{T-\supscr{\tau}{low}}\big) \\
    & \hspace{11em}\ge a 
    \lambda_{\min}\big(\W_{\supscr{\tau}{low}}(t)|_0^{T-\supscr{\tau}{low}}\big),
  \end{align*}
  which concludes the proof.
\end{IEEEproof}

The result of Corollary~\ref{corollary:robust:observability:LTI} is
obtained by the same arguments by using the bound from
Corollary~\ref{corollary:distance:to:Gramian:LTI} instead of that
given in Lemma~\ref{lemma:distance:to:Gramian}.  We are now in
position to prove
Proposition~\ref{prop:state:estimation:output:disturbances}.

\begin{IEEEproof}[Proof of
  Proposition~\ref{prop:state:estimation:output:disturbances}]
%
  %
  %
  To show the result, note that $\bm{\zeta}^i=\Obs_i\xi_{t_{i+\ell}}^i$, 
  or equivalently, $\bm W_i\bm{\zeta}^i=\bm W_i\Obs_i\xi_{t_{i+\ell}}^i$, and 
  consequently $\xi_{t_{i+\ell}}^i=(\bm W_i\Obs_i)^{\dagger}\bm 
  W_i\bm{\zeta}^i$, since $\bm W_i\Obs_i$ is of full rank. Thus, we obtain
  \begin{align}
    \|\widehat{\xi}_{t_{i+\ell}}^i & - \xi_{t_{i+\ell}}^i\| 
    =\| {(\bm W_i\Obs_i)^{\dagger}\bm 
      W_i}(\widehat{\bm{\zeta}}^i-\bm{\zeta}^i)\| \nonumber \\
    & =\|(\bm W_i\Obs_i)^{\dagger}\bm W_i\bm{\delta}^i\| \le \|(\bm W_i 
    \Obs_i)^{\dagger}\|\|\bm W_i\bm{\delta}^i\|. 
    \label{xi:difference:intermediate:bound}
  \end{align} 
  %
  We upper bound the first term on the right hand side of 
  the above inequality as
  \begin{align*}
    \|\bm W_i\bm{\delta}^i\|
    & = \Big(\sum_{l=1}^{\ell}w_l^2(i)\|\delta_l^i\|^2\Big)^\frac{1}{2} \\
    & \le\delta^*\Big(\frac{\tau_1(i)}{2}+\sum_{l=2}^{\ell-1}
    \frac{\tau_{l-1}(i)+\tau_l(i)}{2}+\frac{\tau_{\ell-1}(i)}{2}\Big)^\frac{1}{2}
     \\
    & =\delta^*\sqrt{\tau(i)}\le \delta^*\sqrt{\supscr{\tau}{up}},
  \end{align*}
  where we have use the definition of $w_l(i)$
  from~\eqref{weight:matrix:W} and the bound $\delta^*$ on each
  $\delta_l^i$.  While the second term satisfies
  \begin{align*} 
    \|(\bm W_i\Obs_i)^{\dagger}\|=\sigma_{\max}((\bm W_i
    \Obs_i)^{\dagger})=\frac{1}{\sigma_{\min}(\bm W_i\Obs_i)}, 
  \end{align*}  
  with $\sigma_{\max}$ and $\sigma_{\min}$ denoting the largest and
  smallest nonzero singular value of the corresponding non-degenerate
  matrix, respectively.  The second equality follows from the fact
  that $\bm W_i\Obs_i$ is of full rank (see \cite[Page 435,
  Proposition 4]{PL-MT:85}). Thus, {by taking into
    account~\eqref{eigenvalue:inequality} and that
    $\lambda_{\min}(\Obs_i^{\top}\bm W_i\bm
    W_i\Obs_i)=\sigma_{\min}(\bm W_i\Obs_i)^2$, it follows from
    \eqref{xi:difference:intermediate:bound}}
  that \eqref{xi:difference:bound} is fulfilled.
\end{IEEEproof}

Finally, we give the proof of
Theorem~\ref{thm:numerical:ambiguity:radius:outputs}.

\begin{IEEEproof}[Proof of
  Theorem~\ref{thm:numerical:ambiguity:radius:outputs}]
  The proof is analogous to that of
  Theorem~\ref{thm:numerical:ambiguity:radius}, with the key
  modification being establishment of
  \eqref{distace:between:hatP:barP} with $\bar{\eps}_N$ as in
  \eqref{bar:epsN:output}. By taking into account
  \eqref{numerical:error:bound}, \eqref{xi:difference:bound}, the
  elementary inequality $(a+b)^p\le 2^{p-1}(a^p+b^p)$, which holds for
  any $a,b\ge 0$ and $p\ge 1$ \cite[Lemma 2.4.6]{RBA:72}, and that the
  flow $\Phi$ is linear, as in the proof of
  Theorem~\ref{thm:numerical:ambiguity:radius}, we obtain
  \begin{align*}
    & W_p(\bar P_{\xi_T}^N,\widehat P_{\xi_T}^N)
    \le\Big(\frac{1}{N}\sum_{i=N^{\flat}}^{\bar N}\|\xi_T^i-\bar
    \xi_T^i\|^p\Big)^\frac{1}{p}
    \\
    & =\Big(\frac{1}{N}\sum_{i=N^{\flat}}^{\bar N}
    \|\Phi_{T,t_i^{\ell}}(\xi_{t_i^{\ell}}^i)-\supscr{\Phi}{num}_{T,t_i^{\ell}}
    (\widehat{\xi}_{t_i^{\ell}}^i)\|^p\Big)^\frac{1}{p}
    \\
    & \le \Big(\frac{2^{p-1}}{N}\sum_{i=N^{\flat}}^{\bar N}\big(
    \|\Phi_{T,t_i^{\ell}}(\xi_{t_i^{\ell}}^i-\widehat{\xi}_{t_i^{\ell}}^i)\|^p
    \\
    & \hspace{7em}+\|\Phi_{T,t_i^{\ell}}(\widehat{\xi}_{t_i^{\ell}}^i)
    -\supscr{\Phi}{num}_{T,t_i^{\ell}}(\widehat{\xi}_{t_i^{\ell}}^i)
    \|^p\big)\Big)^\frac{1}{p}
    \\
    & \le \Big(\frac{2^{p-1}}{N}\sum_{i=N^{\flat}}^{\bar N}
    \big((\epsilon^*)^pe^{pL(T-t_i^{\ell})}
    +\fK^p(e^{L(T-t_i^{\ell})}-1)^p\big)\Big)^{\frac{1}{p}}
    \\
    & \le \Big(\frac{2^{p-1}}{N}\sum_{i=1}^N
    \big((\epsilon^*)^pe^{pL\Delta(N-i)} +\fK^p(e^{L\Delta
      (N-i)}-1)^p\big)\Big)^{\frac{1}{p}}
    \\
    & = \Big(\frac{2^{p-1}}{N}\Big[(\epsilon^*)^p\sum_{i=0}^{N-1}
    e^{pL\Delta i} +\fK^p\sum_{i=1}^{N-1}(e^{L\Delta
      i}-1)^p\Big]\Big)^{\frac{1}{p}}
    \\
    & \le \Big(\frac{2^{p-1}}{N}\Big[(\epsilon^*)^p\int_0^N
    e^{pL\Delta s}ds +\fK^p\int_1^N(e^{L\Delta
      s}-1)^pds\Big]\Big)^{\frac{1}{p}} .
  \end{align*}
  By evaluating the first integral in the latter expression we obtain
  the desired result.
\end{IEEEproof}

\bibliography{alias,JC,SM,SMD-add} 
\bibliographystyle{IEEEtranS}

\begin{IEEEbiography}[{\includegraphics[width=1in,height=1.25in,clip,keepaspectratio]{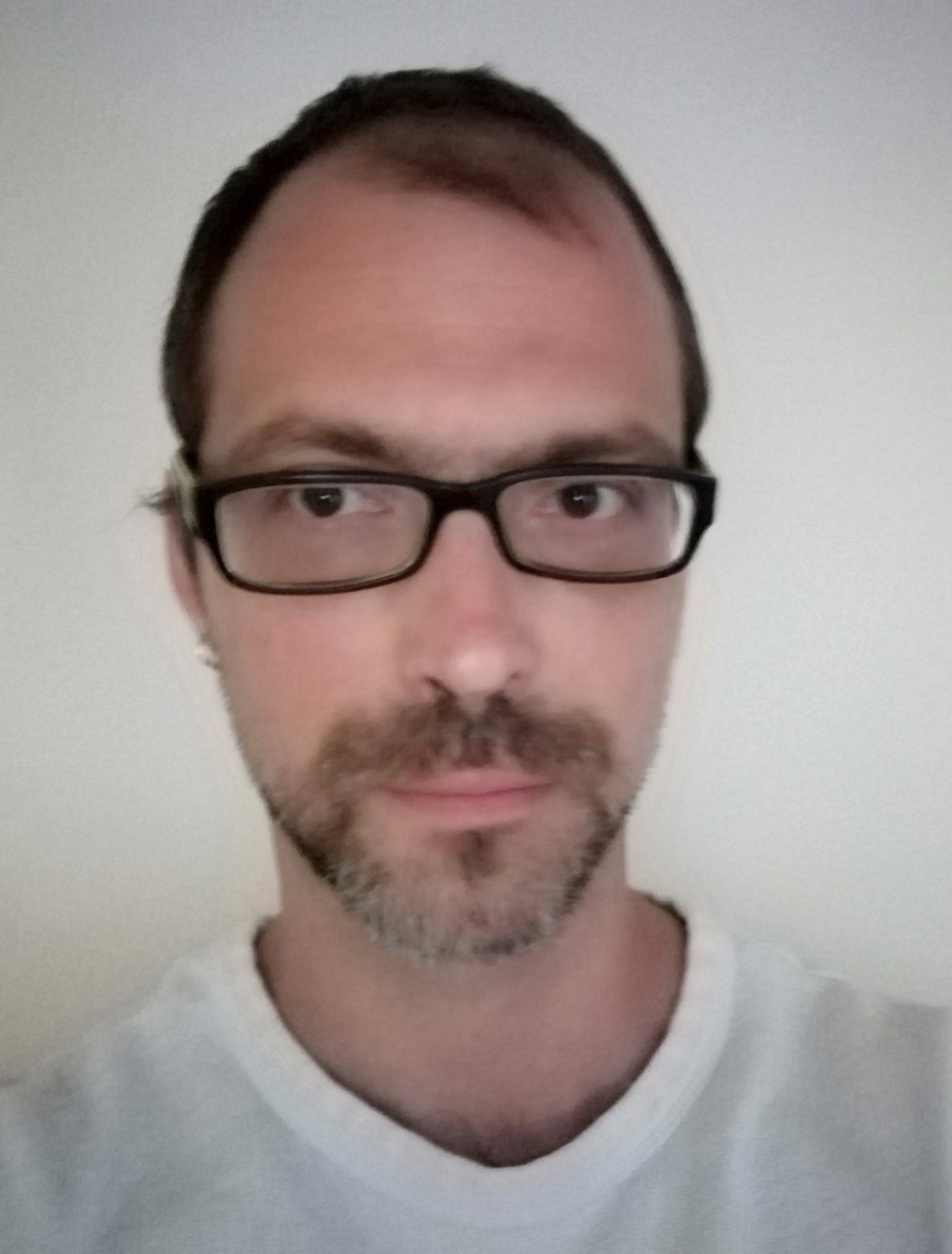}}]{Dimitris
    Boskos} 
  Dimitris Boskos was born in Athens, Greece in 1981. He has received
  the Diploma in Mechanical Engineering from the National Technical
  University of Athens (NTUA), Greece, in 2005, the M.Sc. in Applied
  mathematics from the NTUA in 2008 and the Ph.D. in Applied
  mathematics from the NTUA in 2014. Between August 2014 and August
  2018, he has been a Postdoctoral Researcher at the Department of
  Automatic Control, School of Electrical Engineering, Royal Institute
  of Technology (KTH), Stockholm, Sweden. Since August 2018, he is a
  Postdoctoral Researcher at the Department of Mechanical and
  Aerospace Engineering, University of California, San Diego, CA,
  USA. His research interests include distributionally robust
  optimization, distributed control of multi-agent systems, formal
  verification, and nonlinear observer design.
\end{IEEEbiography}

\begin{IEEEbiography}[{\includegraphics[width=1in,height=1.25in,clip,keepaspectratio]{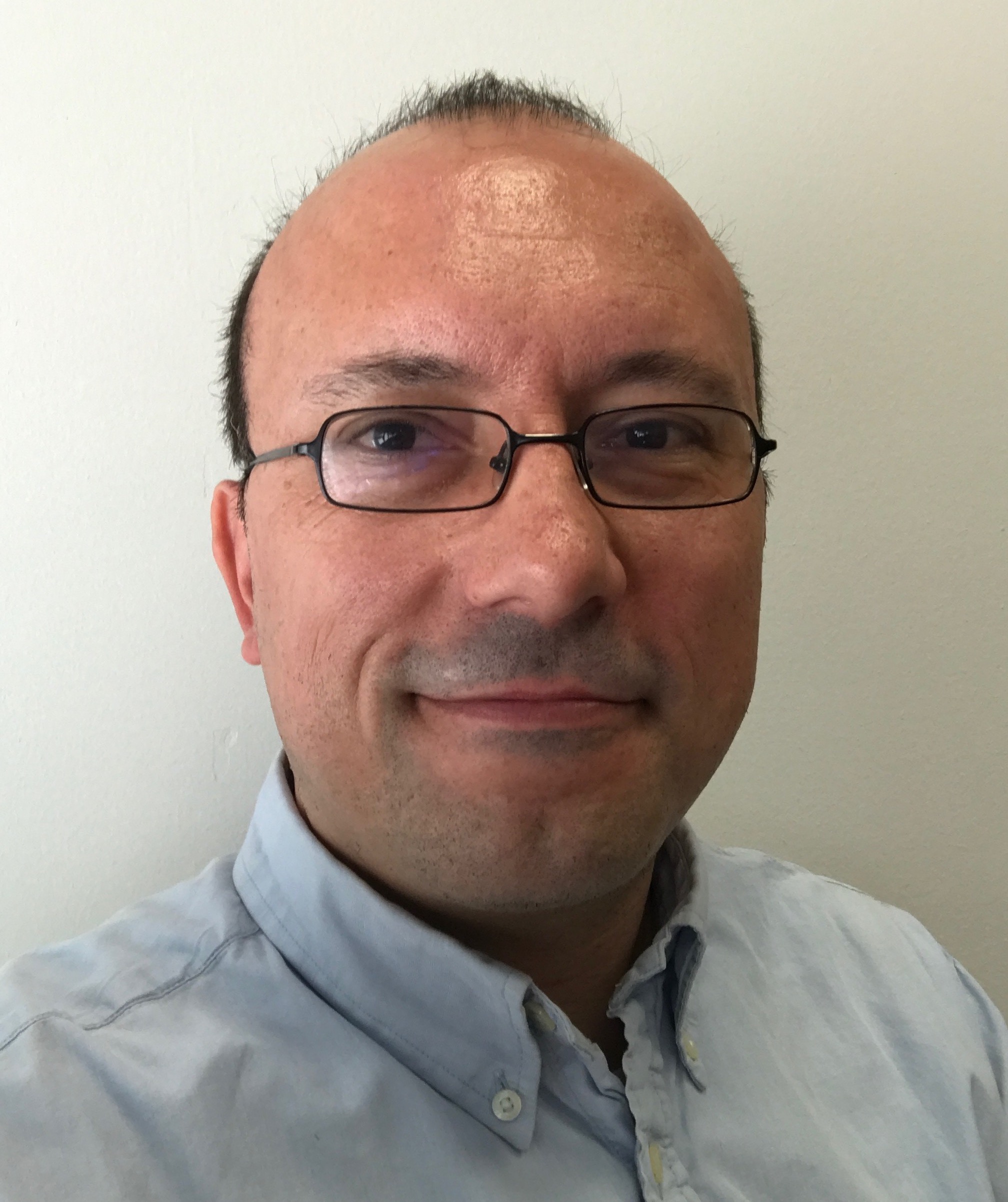}}]{Jorge
    Cort\'{e}s}
  (M'02, SM'06, F'14) received the Licenciatura degree in mathematics
  from Universidad de Zaragoza, Zaragoza, Spain, in 1997, and the
  Ph.D. degree in engineering mathematics from Universidad Carlos III
  de Madrid, Madrid, Spain, in 2001. He held postdoctoral positions
  with the University of Twente, Twente, The Netherlands, and the
  University of Illinois at Urbana-Champaign, Urbana, IL, USA. He was
  an Assistant Professor with the Department of Applied Mathematics
  and Statistics, University of California, Santa Cruz, CA, USA, from
  2004 to 2007. He is currently a Professor in the Department of
  Mechanical and Aerospace Engineering, University of California, San
  Diego, CA, USA. He is the author of Geometric, Control and Numerical
  Aspects of Nonholonomic Systems (Springer-Verlag, 2002) and
  co-author (together with F. Bullo and S.  Mart{\'\i}nez) of
  Distributed Control of Robotic Networks (Princeton University Press,
  2009).  At the IEEE Control Systems Society, he has been a
  Distinguished Lecturer (2010-2014), and is currently its Director of
  Operations and an elected member (2018-2020) of its Board of
  Governors.  His current research interests include distributed
  control and optimization, network science, resource-aware control,
  nonsmooth analysis, reasoning and decision making under uncertainty,
  network neuroscience, and multi-agent coordination in robotic,
  power, and transportation networks.
\end{IEEEbiography}

\begin{IEEEbiography}[{\includegraphics[width=1in,height=1.25in,clip,keepaspectratio]{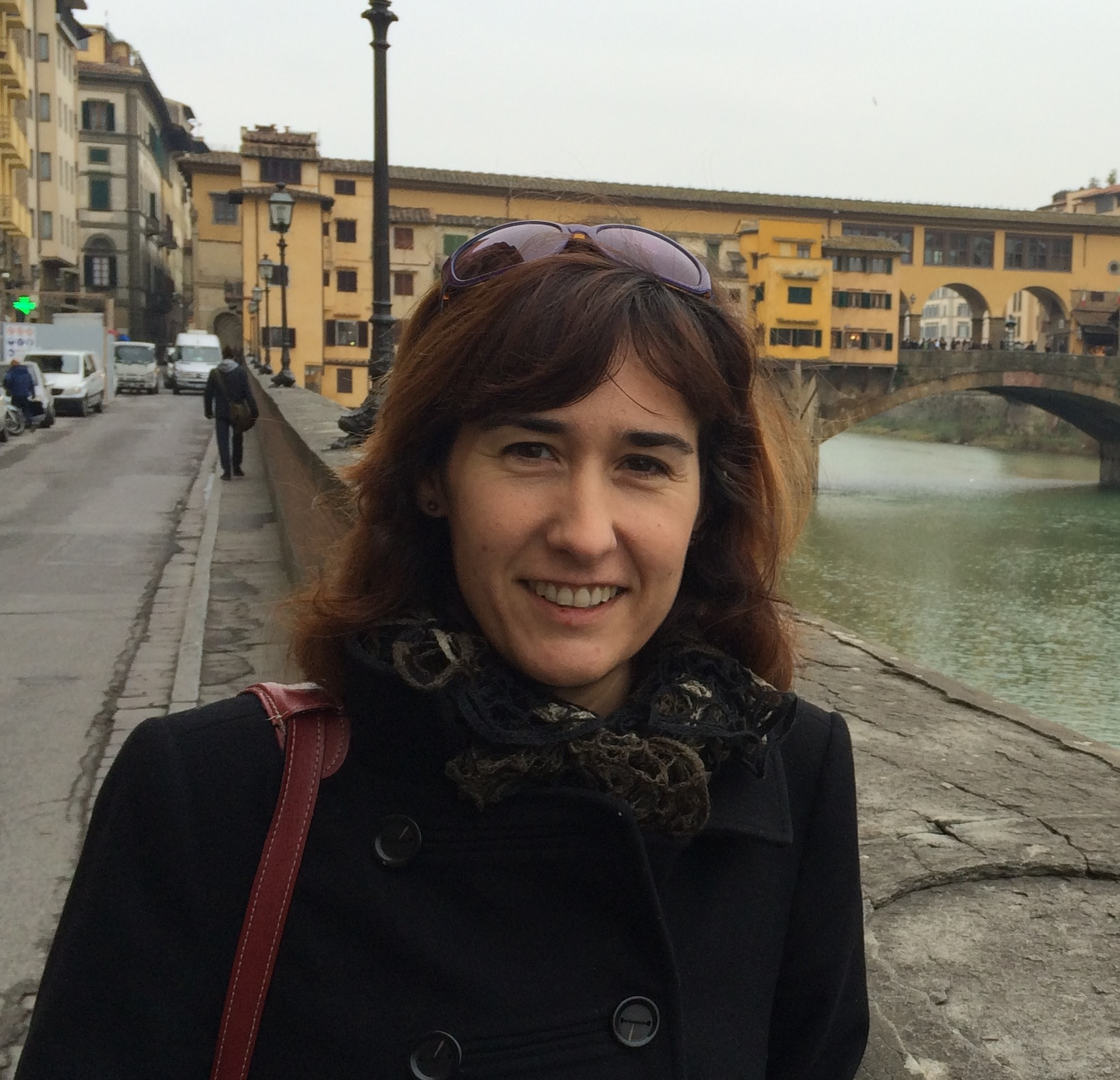}}]{Sonia
    Mart{\'\i}nez} (M'02-SM'07-F'18) is a Professor of Mechanical and
  Aerospace Engineering at the University of California, San Diego,
  CA, USA. She received the Ph.D. degree in Engineering Mathematics
  from the Universidad Carlos III de Madrid, Spain, in May
  2002. Following a year as a Visiting Assistant Professor of Applied
  Mathematics at the Technical University of Catalonia, Spain, she
  obtained a Postdoctoral Fulbright Fellowship and held appointments
  at the Coordinated Science Laboratory of the University of Illinois,
  Urbana-Champaign during 2004, and at the Center for Control,
  Dynamical systems and Computation (CCDC) of the University of
  California, Santa Barbara during 2005.  In a broad sense, her main
  research interests include the control of network systems,
  multi-agent systems, nonlinear control theory, and robotics.  For
  her work on the control of underactuated mechanical systems she
  received the Best Student Paper award at the 2002 IEEE Conference on
  Decision and Control. She was the recipient of a NSF CAREER Award in
  2007. For the paper ``Motion coordination with Distributed
  Information,'' co-authored with Jorge Cort\'es and Francesco Bullo,
  she received the 2008 Control Systems Magazine Outstanding Paper
  Award. She has served on the editorial board of the European Journal
  of Control (2011-2013) and the Journal of Geometric Mechanics
  (2009-present), and currently serves as a Senior Editor of the IEEE
  Transactions on Control of Network Systems.
\end{IEEEbiography}

\end{document}